\documentclass{article} 
\usepackage{enumerate}  
\usepackage{a4}
\usepackage{times}
\usepackage{amsmath}
\usepackage{color}
\usepackage{amssymb}
\usepackage{graphicx}
\usepackage[left=2.75cm,right=2.25cm,top=2.9cm,bottom=2.25cm]{geometry}
\usepackage{booktabs}

\expandafter\chardef\csname pre amssym.def  at\endcsname=\the\catcode`\@
\catcode`\@=11
 
\def\undefine#1{\let#1\undefined} 
\def\newsymbol#1#2#3#4#5{\let\next@\relax
 \ifnum#2=\@ne\let\next@\msafam@\else
 \ifnum#2=\tw@\let\next@\msbfam@\fi\fi
 \mathchardef#1="#3\next@#4#5}
\def\mathhexbox@#1#2#3{\relax
 \ifmmode\mathpalette{}{\m@th\mathchar"#1#2#3}%
 \else\leavevmode\hbox{$\m@th\mathchar"#1#2#3$}\fi}
\def\hexnumber@#1{\ifcase#1 0\or 1\or 2\or 3\or 4\or 5\or 6\or 7\or 8\or
 9\or A\or B\or C\or D\or E\or F\fi}
 \font\tenmsa=msam10
\font\sevenmsa=msam7  
\font\fivemsa=msam5
\newfam\msafam
\textfont\msafam=\tenmsa
\scriptfont\msafam=\sevenmsa
\scriptscriptfont\msafam=\fivemsa 
\edef\msafam@{\hexnumber@\msafam}

\font\tenmsb=msbm10
\font\sevenmsb=msbm7
\font\fivemsb=msbm5
\newfam\msbfam
\textfont\msbfam=\tenmsb
\scriptfont\msbfam=\sevenmsb
\scriptscriptfont\msbfam=\fivemsb
\edef\msbfam@{\hexnumber@\msbfam}

   \font\tengothic=eufm10
   \font\sevengothic=eufm7
   \newfam\gothicfam
   \textfont\gothicfam=\tengothic
   \scriptfont\gothicfam=\sevengothic
   \def\goth#1{{\fam\gothicfam #1}}
   \font\tenmsb=msbm10
   \font\sevenmsb=msbm7
   \newfam\msbfam
   \textfont\msbfam=\tenmsb
   \scriptfont\msbfam=\sevenmsb
   						
\newtheorem{prop}{Proposition}[section] 
\newtheorem{set}[prop]{Setting} 
\newtheorem{rem}[prop]{Remark}
\newtheorem{thm}[prop]{Theorem}
\newtheorem{coro}[prop]{Corollary}

\newtheorem{lemma}[prop]{Lemma}

\newtheorem{ex}[prop]{Example}

\newtheorem{(*)}[prop]{}
\newcommand{\m}{\mbox{${ \goth m}$}}
\newcommand{\B}{\mbox{$\widetilde B$}}

\newcommand{\e}{\mbox{\sf  {$\ell$}}}

\newcommand{{\nat}}{\mbox{${\rm I\hspace{-.6 mm}N}$}}
\newcommand{{\natsmall}}{\mbox{${\rm\scriptstyle I\hspace{-.6 mm}N}$}}

\newcommand{\I}{\mbox{$\  \Longrightarrow \ $}}
\newcommand{\II}{\mbox{$\  \Longleftrightarrow \ $}}

\newcommand{\enu} {\begin{enumerate}} 
\newcommand{\enua} {\begin{enumerate}[ $(a)$]} 
\newcommand{\enui} {\begin{enumerate}[ $(i)$]} 
\newcommand{\denu} {\end{enumerate}}

\font\tengothic=eufm10 
\font\sevengothic=eufm7
\newfam\gothicfam 
\textfont\gothicfam=\tengothic
\scriptfont\gothicfam=\sevengothic
\def\goth#1{{\fam\gothicfam #1}}
\begin{document} 
 \begin{center}

\renewcommand{\thefootnote}{\fnsymbol{footnote}}
\addtocounter{footnote}{0}

{\LARGE   On semigroup rings with  decreasing Hilbert function.}
\end{center}

\begin{center}
{\large
  Anna Oneto$^1$,   Grazia Tamone$^2$. 
}  \\[2mm]
{\small $ $Dima - University of Genova,
via Dodecaneso 35, I-16146 Genova,  Italy.}\\
{\small E-mail:  $^1$oneto@dime.unige.it, \ $^2$tamone@dima.unige.it }
\end{center}

\vspace*{-3mm}

\noindent \rule[0pt]{\textwidth}{1pt}

\vspace*{-2mm}

{\baselineskip8pt

\begin{abstract}

{\baselineskip5pt

\noindent In this paper  we study the Hilbert function $H_{\!R}$ of  one-dimensional semigroup rings $R=k[[S]]$.  For some classes of semigroups, by means of the notion of {\it support} of the elements in $S$,  we give  conditions on the generators of $S$ in order to have decreasing $H_{\!R}$. When the embedding dimension $v$ and the multiplicity $e$ verify $v+3\leq e\leq v+4$,   the decrease of $H_{\!R}$ gives   explicit description  of  the  Ap\'ery  set of $S$.  In particular for $e=v+3$, we classify the semigroups with $e=13$ and $H_{\!R}$ decreasing,  further we show that $H_{\!R}$ is non-decreasing if $e<12$.  Finally we deduce that   $H_{\!R}$ is non-decreasing  for every Gorenstein semigroup ring   with $ e\leq v+4$.  }
\end{abstract}
}

\vspace*{-2mm}

\noindent \rule[0pt]{\textwidth}{1pt}

{\baselineskip6pt
\vskip4pt

\noindent {\footnotesize\it Keywords~}{\scriptsize\sf :} {\scriptsize\sf
 Numerical semigroup,   Monomial curve, Hilbert function, Ap\'ery  set. }
\smallskip

\noindent {\footnotesize\it   Mathematics Subject Classification~}: {\scriptsize\sf Primary~:
14H20~; Secondary~: 13H10.}

}

\vspace*{-5mm}

\renewcommand{\thefootnote}{\arabic{footnote}}
\addtocounter{footnote}{0}

\setcounter{section}{-1}

\section{Introduction.}
\indent Given a local noetherian ring $(R,\m,k)$ and the associated graded ring $G=\oplus_{n\geq 0}\big( \m^n
/\m^{n+1}\big)$, a classical hard topic in commutative algebra is the   study of the Hilbert function $H_R$ of $G$, defined as $H_R(n)=dim_k\big( \m^n
/\m^{n+1}\big)$: when $R$ is the local ring of a $k$-scheme $X$  at a point $P $, $H_R$ gives important geometric information.
  If $depth(G)$ is large enough, this function   can be computed by means of the Hilbert function of a lower dimensional ring, but in general $G$ is not Cohen-Macaulay, even if $R$ has this property.\par 
	For a   Cohen-Macaulay one-dimensional local ring $R$, it is well known that $H_R$ is a non decreasing function when  $G$ is  Cohen-Macaulay, but we can have $depth(G)=0$ and   in this case $H_R$ can be decreasing, i.e. $H_R(n)<H_R(n-1)$ for some $n$   (see, for example    \cite{gr},\, \cite{o},\, \cite{r}).  This fact cannot happen if $R$ verifies  either $v\leq 3,$ \, or $\,v\leq e\leq v+2 $, where   $e$ and   $v$ denote respectively the  multiplicity and the embedding dimension of $R$    (\,see  \cite{e1},\, \cite{e2}, \,\cite{rv}\,). \, 
If $R=k[[S]]$  is a semigroup ring,   many  authors proved that  $H_R$ is non-decreasing in several cases: \, 
$\bullet$  $S$ is generated by an almost arithmetic sequence (if the sequence is arithmetic, then $G$ is Cohen-Macaulay) \cite{t}, \cite{mt1}\,\, 
$\bullet$ $S$ belongs to particular subclasses of four-generated  semigroups, which are symmetric\,\,\cite{am},  or which have  Buchsbaum tangent cone \,\cite{cjz} \, \,
$\bullet$ $S$ is  balanced \,\,  \cite{pt},\, \cite{cjz}\,\,
$\bullet$ $S$ is obtained by techniques of gluing numerical semigroups \,\, \cite{ams},\,  \cite{jz}\,\,
$\bullet$   $S$ satisfies certain conditions on particular subsets of $S$ (see below) \cite[Theorem 2.3, Corollary 2.4, Corollary 2.11]{ddm}.   \, 
If  $e\geq v+3$, \, the function $H_R$ can be decreasing, as shown in several examples: the first one (with $e=v+3$) is   in  \cite{mt2}  (here recalled in  Example \ref{tammol}).
  When $G$ is not Cohen-Macaulay, a useful method to describe $H_R$   can be found in some recent papers (see \cite{pt}, \cite{cjz},\cite {ddm}): it is based on     the study of   certain subsets of $S$, called $D_k $ and $C_k, \, (k\in {\mathbb N})$. \par
  The aim of this paper is the study  of  semigroup rings   $R=k[[S]]$ with $H_R$ decreasing.  To this goal we introduce and use the notion of   {\it support}  of the elements in $S$ (\ref{def2}.4); by means of this tool   we first develop a technical  analysis of the  subsets $D_k,\, C_k$   in Section 2. \, Through this machinery, under suitable assumptions on the Ap\'ery  set of $S$, we prove (Section 3) necessary conditions on $S$ in order to have decreasing Hilbert function,  see  (\ref{apj1}), (\ref{ap24}).\\  
 In Section 4 we apply  these results to the semigroups    with $v\in\{e-3,e-4\}$. \par For $v=e-3$ we show that the decrease of $H_R$ is characterised by a particular structure of the sets $C_2,\, D_2, \, C_3$ and that    $H_R$ does not decrease for $e\leq 12$, see  (\ref{e-3}), (\ref{e13});  in particular, for $e=13$, we identify precisely   the semigroups with  $H_R$ decreasing, see    (\ref{exe13}) and examples (\ref{es13}).\\
 In case  $v=e-4$ we obtain analogous informations on the structure of $C_2, C_3, D_2, D_3$, see  (\ref{40}) and (\ref{31}).\par 
 Such methods allow to construct various examples of semigroup rings with  decreasing $H_R$, see, for example   (\ref{es2}), (\ref{es5}) where  $e-7\leq v\leq e-3$;   in particular   example (\ref{es5}.1)   describes    a semigroup whose Hilbert function decreases at two different levels.
 \,\, The examples have been performed by using the program CoCoA together with FreeMat and Excel.  \par
 As a consequence of some of the above facts, one can see that  the semigroups $S$ with $| C_2|=3  $ and $| C_3\cap ${\it Ap\'ery} $ \,set |\leq 1$  cannot be symmetric.    It follows   that  every Gorenstein ring  $k[[S]]$ with $\,v\geq e-4\,$   has   non-decreasing Hilbert function (\ref{corogor}). This result  is a partial answer to the    conjecture  settled by M.E. Rossi \cite[Problem 4.9]{ro} that  a   Gorenstein $1$-dimensional local ring has non-decreasing Hilbert function. 
    
 \section{Preliminaries.}
 
   We breafly recall the definition of the Hilbert function for local rings. Let  $(R,\m,k)$ be a noetherian local $d$-dimensional ring,  the {\it associated graded ring} of $R$ with respect to $\m$ is $G:= \bigoplus_{n\geq 0} \m^n/\m^{n+1} $:
the {\it Hilbert function} $H_R: G\longrightarrow {\mathbb N}$ of $R$ is   defined by    $H_R(n)= dim_k(\m^n/\m^{n+1})$.\,\,This function   is called \, {\it non-decreasing}\, if $ H_{R}(n-1)\leq H_R(n)$ for each $n\in{\mathbb N}$ and  \,  {\it decreasing}\,  if there exists $\e\in{\mathbb N}$ such that $H_R(\e-1)>H_{R}(\e )$, \,  we say $H_R$ {\it decreasing at level} $\e$.  

   \begin{(*)}\label{CM} 
{\rm   Let $R$ be    a one-dimensional   Cohen-Macaulay local ring 
 and assume  $k=R/\m$ infinite. Then there exists a {\it superficial element} $x\in \m$, of degree 1,   (i.e. such that $x\m^n=\m^{n+1} $   for $n>>0$).\\[2mm]
 It is well-known that
 \begin{itemize}
 \item  $G$ is Cohen-Macaulay $ \II $ the image $x^*$ of $x$ in $G$ is a  non-zero divisor. 
\item If $G$ is Cohen-Macaulay, then  $H_R$ is non-decreasing. 
  \end{itemize}} 
\end{(*)}
  We begin by setting the notation of the paper and by recalling some known   useful facts. 
 \begin{set}\label{def1} 
{\rm In this paper $R$ denotes a 1-dimensional numerical semigroup ring, i.e.  $R=k[[S]]$, where $k$ is an infinite field and  $S=\{\sum a_in_i, \,a_i,n_i\in{\mathbb N}\}  $ is a {\it numerical semigroup } of {\it multiplicity e} and {\it embedding dimension} $v$  minimally generated by $\{e:=n_1,n_2,\dots ,n_{v }\}$, with $ 0<n_ 1<\cdots <n_{v },\,\, GCD(n_1,\dots ,n_v)=1$. Then $R$ is the completion of the local ring   $\,k[x_1,\dots,x_{v }]_{ (x_1,\dots ,x_{v })}$ of the monomial curve  ${\cal C}$ parametrized by $x_i:=t^{n_i}\, \, (1\leq i\leq v )$. The maximal ideal of $R$ is     $\m=(t^{n_1},\cdots ,t^{n_v})$ and  $x_1 =t^e$ is a superficial element of degree 1.\ Let ${\rm v} :k((t))\longrightarrow {\mathbb Z}\cup\{\infty\}$ denote the usual valuation. } 
\enu 
\item  {\rm $M:=S\setminus\{0\}={\rm v}(\m),\quad hM ={\rm v}(\m^h) ,\,\, $ for each $h  \geq 1$  \,  \, and the Hilbert function $H_R$ verifies
  \item[]  $H_R(0)=1$,\quad $H_R(h)=|   hM\setminus (h+1)M  | ,\,$ for each $h\geq 1 $.}
\item{\rm  Let $g\in S$,   the  \, {\it   order} \, of $g$ is defined as\,
$ord(g):=max\{ h\,|  \, g\in hM \}$.}
\item {\rm {\em Ap\,=\,Ap\'ery(S)}$ \,:=\{ s\in S\,|\, s-e\notin S\}  $ is {\it the Ap\'ery  set of S}  with respect to the multiplicity $e$ , $| Ap|=e $ \, and $\, e+f\,$   is the greatest element  in $Ap$, where $f:=max\{x\in\nat\setminus S\}$ is    the {\it Frobenius number} of $S$. \\
Let $d:=max\{ord(\sigma) \,|\,\, \sigma\in Ap\}$.\\ Denote by
  $Ap_k:=\{ s\in Ap \, |\,ord(s)=k\},\,\, k\in[1,d]$   the subset of the elements of order $k$ in {\it Ap\'ery}.}
  
\item {\rm  Let $R':=R/t^eR$, the Hilbert function of $R'$  is \, $H_{R'}=[1,a_1,\dots,a_d]$ with $a_k=| Ap_k| $ \, for each $k\in[1,d]$, \, see, for example , {\rm \cite[Lemma 1.3] {pt}}. }
 \item {\rm A semigroup $S$ is called {\it symmetric} if for each $s\in S$  \quad   $s\in Ap\II e+f-s\in Ap$. }
  \denu   
  \end{set} 
  By (\ref{def1}.1)  and (\ref{CM}), if $ord(s+e)=ord(s)+1$ for each $s\in S$, then $G$ is Cohen-Macaulay and $H_R$ is non-decreasing. Therefore in order to focus on the possible decreasing Hilbert functions, it is useful to define the following subsets $D_k,C_k\subseteq S$,  we also introduce the notion of  {\it support }  for a better understanding of these sets.
  \begin{set}\label{def2} 
\enu 
 
   \item {\rm  $D_k:=\{s\in S\,|\, ord(s)=k-1\, and\,  ord(s+e)>k\}.\, \big(D_1=  D_k=\emptyset,\, \forall k\geq r  $ {\rm \cite[Lemma 1.5.2] {pt}}$\big)$.
 \item[] $D_k^t:=\{s\in D_k$ {\it  such  that\, } $ord(s+e)=t\}$ \, and let 
  \,\, $k_0 :=  min\{ k\,\,{\it such\, that} \,\,D_k\neq\emptyset\}$.}
\item  {\rm $C_k:=\{s\in S\,|\, ord(s)=k\,\,\,and\,\,\, s-e\notin (k-1)M\}$, \,\, i.e.,
  $C_k= Ap_k \, \bigcup \, \big\{  \cup_h (D_h^k+e), \,\, with \,\,2\leq h\leq k-1\big\}.$
\item[]   $      C_1=\{n_2,\cdots ,n_v\},\,\, C_2= Ap_2,\,\,  C_3=(D_2^3+e)\cup Ap _3 $ \,   and \,\,
 $ C_k =\emptyset,  \,\,\forall k\geq r+1  $\,
 {\rm \cite[Lemma 1.8.1] {pt},  \,\cite{cjz}}$ $. }
 
  \item {\rm  A  {\it  maximal representation}  of $s\in S$ is     any expression $  s=\sum_{j=1 }^v a_jn_j , a_i\in\nat$,   with $  \sum_{j=1 }^v a_j=ord(s)$ \item[] and in this case  
   we define  {\it  support }  of    $s $   as  \, $Supp(s):=  \{n_i\in Ap_1\,|\, a_i\neq 0\}. \,  $ Recall: \\ $Supp (s)$   depends on the  choice of a maximal representation  of $s  $. 
 Further for $s\in S$, we define:\\   $|Supp(s)|:= max\big\{ |Supp\left(\sum_i a_in_i\right)|,  \,{\it such\,that} \,  s=\sum_i a_in_i $   {\it   is  a   maximal representation of s}$\big\}$.}
 
      \item {\rm For a subset $H\subseteq S$,  \,  $Supp (H):= \bigcup \big\{  Supp (s_i ), s_i\in H\big\}$.}

      \item 
   {\rm   We   call   {\it   induced by} $s=\sum_{j=1 }^v a_jn_j$      (maximal representation) an element \par$\, s'=\sum_{j=0 }^v b_j n_j$, with \,\, $0\leq b_j\leq a_j$. \quad    Recall:     $ord(s')=\sum b_j$\, \, by    \cite[Lemma 1.11] {pt} }.  

 \denu 
\end{set}
The following two propositions are   crucial in the sequel. 
\begin{prop}\label {PT0}   
Let $S$  be as in Setting \ref{def1} and let $s  \in C_k, \, k\geq 2$. Then 
\enu
\item For each  $s=\sum 
_{j\geq 1} a_jn_j\in C_k \,\,  ($maximal representation, with $ \sum 
_{j\geq 1} a_j=k )$, we have:
\begin{enumerate} 
\item  Every element $s'=\sum 
_{j\geq 1} b_jn_j   $  induced by $s$, with $\sum 
_{j\geq 1} b_j=h$,\, belongs to $C_h$. 
\item   $Supp \, (C_k)\subseteq Supp (C_{k-1})\subseteq...\subseteq Supp(C_2)$.
  \denu 
     \item If $ s=g+e$ with $g\in D_h, \, (h\leq k-1)$,     any maximal representation    $s= \sum 
_{j\geq 1} a_jn_j$ has $a_1=0$.

\denu
\end{prop}
Proof. ({\it 1.a}). See {\rm \cite[Lemma 1.11] {pt}}.\\
({\it 2}). If $a_1>0$, then $g =(a_1-1)e+\sum 
_{j\geq 2} a_jn_j$ has $ord=k-1$, by ({\it 1.a}), contradiction.\quad$\diamond$

\begin{prop}\label {pt}    
Let $S$  be as in Setting \ref{def1},   let $k_0$ be as in (\ref{def2}.1) and let \, $k\geq 2$.\enu
 \item $H_R(k)-H_R(k-1)=| C_k|- |D_k|$ \,   {\rm \cite[Proposition 1.9.3]{pt}, \cite[Remark 4.1]{cjz}}
 \item $G$ is Cohen Macaulay $\II D_k=\emptyset$  \, for each $k\geq 2$.\,  {\rm\cite[Theorem 1.6] {pt}}.
 \item If $| D_k|\leq k+1$ for every $k\geq 2$, then $H_R$ is non-decreasing \,  {\rm \cite[Theorem 3.3, Corollary 3.4]{ddm}}.

\item   If $ | D_k|  \geq k+1$, then $\ | C_h|\geq h+1 $ \, for all \, $  h\in[2,k]$ \,  {\rm \cite[proof of Proposition 3.9]{ddm}}.
 \item  In particular $($ recall: $k_0=min\{k\,|\, D_k\neq \emptyset\}\,)$ {\rm \cite[ Corollary 3.11]{ddm}}$:$ \begin{enumerate}  
 \item If $H_R$ is decreasing at level $k$, then $| D_k|\geq max\{1+ | C_k|,\, k+2\},$ 
 \item If $H_R$ is decreasing, then $| Ap_{k_0}|\geq k_0+1,$ 
  \item If $H_R$ is decreasing, then $|Ap_{2}|\geq 3$.
\denu
 
\denu
\end{prop}
We show a semigroup $S$ with $H_R$ decreasing,  which  is the {\it first example} in the sense that $e=v+3=13,\,\,| Ap_2|=3$  are minimal in order to have decreasing Hilbert function, see (\ref{e-3}),  (\ref{e13}).
\begin{ex}\label{tammol}   
{\rm \cite[Section 2]{mt2}} \quad  Let \quad$S=<13,19,24, 44, 49,54,55,59,60,66>$, with $e=13,\,v=10,$\\
{\rm $Ap =  \{\begin{array}{cccccccccccccccc}
  0& 19,& 24,&  {38},&  {43},&  {44},&  {48},&  {49},&  {54},&  {55},&  {59},&  {60},&  {66}\end{array}\}$.\\
\centerline{$ \begin{array}{rccccccccccccccccc}
 M\setminus2M&  13& { 19} & { 24}&   &   & \bf{44}&   & \bf{49}& \bf{54}&  {55}&  \bf{59}&  {60}&  {66}& dim&=& 10\\
2M\setminus3M& 26& 32& 37& \bf{38}& \bf{43}&   & \bf{48}&   &   & 68&   & 73& 79& dim&=& 9\\
3M\setminus4M& 39& 45& 50& 51& 56& 57& 61& 62& 67&   & 72&   & 92& dim&=& 11\\
4M\setminus5M& 52& 58& 63& 64& 69& 70& 74& 75& 80& 81& 85   & 86& & dim&=& 12 \\
5M\setminus6M& 65& 71& 76& 77& 82& 83& 87& 88& 93& 94& 98& 99& 105& dim&=& 13\\
\end{array} $}\\[1mm]
 Then $H_R= [1,10,9,11,12,13 \rightarrow]$, \,\, $H_{R'}= [1,9,3]$. Further   $Supp(D_2+e)=Supp\,C_2=\{19,24\}$}: \\[1mm]
$\begin{array}{lllllllllccc} 
 {\bf D_2} &=&{\bf \{44,49,54,59\}} &{\bf C_2 }={\bf\{38,43,48\}}=\{19\cdot2,19+24,24\cdot2\} =Ap_2\\[1mm]
  D_2+e &=&\{57 ,62 ,67 ,72 \}&57 = 3\cdot19,\quad 62= 2\cdot19+24,\quad 67=  19+2\cdot24,\quad 72= 3\cdot24 \\[1mm]

 {  D_3} &=&\{68,73\}  \, 
 & {  C_3} 
 = \{57,62,67,72\} =D_2+e\\[1mm]
 
D_3+e&=&\{81,86\}& 
81= 3\cdot 19 +24,\quad 86= 2\cdot 19 +2\cdot 24\\[1mm]

 {   D_4} &=&\{92\} 
 &{   C_4 }
 =\{81,86\} =D_3+e\\[1mm]
D_4+e&=&\{105\} &105=3\cdot19+2\cdot24\\[1mm]

 {  D_5}&=&\{\emptyset\}  & {  C_5 }=\{105\}=D_4+e.\\
\end{array}$    
\end{ex} 
  \section{Technical analysis of $\, C_k \,$ and $\,D_k$     via    supports and   {\it Ap\'ery } subsets.}
\begin{lemma}\label{suppq} 
    Let   $x=\sum_{i=2}^v  a_{i }n_{i }\in C_k$, with $a_{i }\geq 1$ , for each $i ,\,\sum_{i=2}^v   a_{i } =ord(x)=k$ and $| Supp(x)|=q$.\\[2mm] Let \quad $2\leq h<  k$,\quad then      \quad$\left[\begin{array}{llllll}
   (a)&| C_h|\geq hp+1\geq q, &  if & q \geq h+1,\, p\geq 1\\
 
 (b)& | C_h|\geq q,&  if & q\leq  h.   
  \end{array}
 \right.$ 
 \end{lemma}
 Proof.  First recall that, by  (\ref {PT0}.1$a$), every  element   of order $h$  induced by the given maximal representation  of $x$, belongs to $C_h$. 
    We denote for simplicity   $Supp(x)= \{m_1<m_2<\cdots .<m_q\}$,   distinct minimal generators, with $m_i\neq e $ by (\ref{PT0}.2).\\[2mm]
    $(a)$.   \ If $q\geq h+1$.\\ Then we can construct the $ (h+1)+h(q-h-1)$ distinct  induced elements in $C_h$: \\[2mm]
  \centerline{$\left\{\sigma_{\eta,i}=(\sum_{j=\eta}^{\eta+h}  m_{j})-m_i,\quad  \eta=1,\cdots ,q-h,\, \left[\begin{array}{llll}
  i= 1,\cdots ,h+1, & {\rm for}   &\eta=1 ,\\
   i= \eta,\cdots ,h+\eta-1,&{\rm for} & 2\leq\eta\leq q-h, 
   \end{array}
 \right. \,\right\}$.} \\[2mm]
$(b)$.  If $1\leq q\leq h  $,    there exists    $(a'_{1 },\cdots ,a'_{ q}) $ such that $a_i\geq a'_{i }\geq 1,\,\forall\,i=1,\cdots ,q$,  $\, \sum_i a'_{i }=h+1$.   Then   $C_h$ contains   the $q $  distinct elements \,
{$\sigma_j =\sum_1^{q} a'_{i }m_{i }-m_j,\,\, j=1,\cdots ,q$.$\quad\diamond$ }

 \begin{prop}\label{CnonAp}  
  Let $k_0=min \{k\in {\mathbb N} \,|\,D_k\neq \emptyset\}, \, d =max\{ord(\sigma), \,\, \sigma\in Ap\}$. 
   \enu
\item Let  $g\in D_{k },   \, g+e=\sum 
_{j} a_jn_j,$   with $\, \sum 
_{j} a_j =k+p, \,\, p\geq 1$ $($maximal representation$)$: 
\enua
\item Let    $y\in C_h, \, h<k+p,$ be induced by $g+e;$ if  \, $  h\leq  max\{ p+1,k_0\}$, \, then  \,\! $y$  \,belongs to $Ap$. 
 \item  $p\leq d-1$.
 \item    $|Supp(g+e)|\leq |Ap_{p+1}|$.
 \denu
\item   If \, $Ap_3=\emptyset$, \, then \,\,  $D_k+e=C_{k+1}$\, \, for each $k\geq 2$.
 
\denu
\end{prop}
 Proof.  {\it 1.$(a)$}\, Let  $  \, h\leq k_0 $ and let    $z\in C_h\setminus Ap; \,$  then $z-e\in S,\,  $ with \, $ ord(z-e)\leq h-2$, hence $z-e\in D_{r},\, r<k_0$,  impossible.  
\,\, Further, if $y=\sum b_jn_j$ and $y\notin Ap,$ with $ \, h=\sum b_j \leq p+1$, then $y= e+\sigma $,  with \, $0<\sigma \in S$.\,\, Then:   {$ g= \sigma+\sum (a_j-b_j)n_j\I ord (g)\geq 1+k+p-(p+1)=k,$ \quad  contradiction.}
\\[2mm]
{\it 1.$(b)$}\,    Let $  g+e=\sum_j a_jn_j $ (maximal representation), if $ \sum_j a_j\geq k+d$, then $d+1\leq \sum_j a_j -k+1$ and so every induced element $\in C_{d+1}$ belongs to $Ap$ by (1), impossible, since $(d+1)M\subseteq M+e$; hence $ord(g+e)=k+1$.
\\[2mm]
{\it 1.$(c)$} Let $| Supp(g+e)|=q$ and let  $h:=p+1<p+k$, by (\ref{suppq}) there are at least $q$ elements in $C_h$. These elements are in $Ap$, by 1.($a$). Hence $q\leq | Ap_{p+1}|$.
\\[2mm]
{\it 2}. \, If there exists $g\in D_k$ such that $ord(g+e)\geq k+2$, then $p\geq 2;$  by 1.($a$) we would have $y\in Ap_3$ for every $y\in C_3$ induced by $g+e$.  \quad $\diamond$
 
  \begin{rem}\label{short}  
\enu   \item {\rm (\ref{CnonAp}.1 $a$)
 cannot be improved, for example in (\ref{tammol}) if $s= 92: \quad ord(s)=3$,    \,\, $92+e=105=3\cdot 19+ 2\cdot 24\I 92\in D_4$, \,\,$p=1$ \,\, and $ 2\cdot 19+   24=62=49+e\in C_3\setminus  Ap_3$ \,   (here $3=p+2$).} 
   \item {\rm Let $x_1,x_2\in D_k, \, x_1\neq x_2$ such that $Supp(x_1+e)=Supp(x_2+e)$,   \,\,\ \
$x_1+e=\sum \alpha_in_i,$\\ $  x_2+e=\sum \beta_in_i$.  Then there exist \,\, $i,j$ \,\, such that $\alpha_i>\beta_i,\quad \alpha_j<\beta_j$.}\par 
  {\rm Proof.  If  $x_1\neq x_2$ and  $\alpha_i\geq \beta_i$ for each $i$, then $  x_1+e=x_2+e+\sigma$, $ord(\sigma)\geq 1$ , hence $ord(x_1)>ord(x_2)$, impossible. }
\item Let $x\in D_k,\,\, ord(x+e)\geq k+2$. 
 If $y\in D_k$ and $ord(y+e)=h\leq k+1$. Then $y+e$ cannot be induced by $x+e$.\par{\rm 
 Proof.  If $x+e=y+e+s$, with $ord(s)\geq 1$ then $x=y+s$, $ord(x)>ord(y)$, impossible, since by assumption  $ord(x)=ord(y)=k-1$.}
\denu  
\end{rem}
Given the  sets $C_h,C_k$, with $h<k $, we   estimate lower bounds for the cardinality of     $  C_h$, by enumerating the  elements induced by $C_k$. We first consider the elements induced by the subset $\{x\in C_k $ {\it such that} $| Supp(x )|\leq 2\}$.
 
 \begin{lemma}\label{supp2}  
  For  $x_r \in C_k,$  $ (k\geq 3)$,  \,let \, $x_r = a_r n_i+b_r  n_j,\, a_r+b_r=k $  and let  \,      $r\in[2, k-1]$.    \enu 
  \item Let $x_1 \in C_k;$   the number  of distinct elements of $C_h$ induced by $x_1$ is  \,
  
   \item[]   $\beta_{ 1} =1+  min\{ a_1,b_1,h, k-h\} =   \left[\begin{array}{llllllll}   1+min\{b_1,h\},&if& a_1\geq h\, &\\ 
  1+a_1&if& a_1<h \leq b_1&\\ 
    1+k-h &if& a_1,b_1<h &\\ 
   \end{array}\right.$ \\
  with  $  \beta_{ 1}\geq 1 \,\, $ and $  \beta_{ 1}=1\II  \, a_1b_1=0$.    \item   Let $x_1,x_2$ be distinct elements in $C_k$ such that \,   $  Supp(x_1 )\cup  Supp(x_2 ) =\{n_i,n_j\}$.\\ We can assume that  
 $\left\{\begin{array}{llll}x_1 = a_1 n_i+b_1  n_j,   \\ 
 x_2 = a_2 n_i+b_2 n_j,\\  a_i+b_i= k,\quad 0\leq a_1<a_2 , \, \I\, 0\leq  b_2<b_1\,.  \\
\end{array}\right.$ \\ [2mm]
The number of distinct elements of $C_h$ induced by $\{x_1,x_2\}$ is \\[2mm]
$\beta_{ 2} = \left[\begin{array}{lllrrllll}
 (a)& 1+  min\{b_1,h\}  \geq 2\,\, (=2\II   b_1=1 , b_2=0)  & if&   h\leq a_1<a_2 \\[1mm]
 
(b) &  2+ min\{ a_1, k-h \}+min\{  h\!-\!a_1 -1, b_2\}\geq 3& if & 0\leq a_1<h\leq a_2\\
(c)&  2+ min\{ a_1, k-h \}+min\{ a_2\!-\!a_1\!-\!1 , k-h\} \geq 2 & if & 0\leq a_1<a_2<h \\
& \beta_{ 2}  =2\II a_1=0,a_2=1
 \end{array}\right.$ 
 \item   If $\,\,   C_k \supseteq \{x_1,x_2,x_3\}$, \,\, let   $\left\{\begin{array}{llll}x_1 = a_1 n_i+b_1   n_j   \\ 
 x_2 = a_2 n_i+b_2 n_j,&    \\
 x_3 = a_3 n_i+b_3 n_j,& \\ a_i+b_i= k, \quad 0\leq a_1<a_2<a_3, \quad b_1>b_2>b_3\geq 0    \end{array}\right.,$  \\[2mm]
 The number of distinct elements of $C_h$ induced by $\{x_1,x_2,x_3\}$  is \, 
  
 $ \beta_{ 3}=  \left[\begin{array}{lrllllll}
  \!\!(a) \,\,\,1+ min\{b_1,h\} \geq 3 & if&   h\leq a_1< a_2 <a_3   \\[2mm]
 \!\! (b) \,\,\,  2+min\{ a_1, k-h \}+min\{h\!-\!a_1\!-\!1 , b_2\! \} \geq 3  & if & a_1<h\leq a_2<a_3  \\[2mm]
\hspace{0.5cm}further\,\, \beta_3\geq 4, \,\,  except\, \,the\,  \,cases: \\ \hspace{4cm}     (i)\,\,\,\,   h=2,\,\, a_1\in\{0,1\}   \\
 
\hspace{4cm} (ii)\,\,\,\,   k=h+1,\,\, a_1=h\!-\!1,\,   b_1=2 \ \\
 \hspace{4cm} (iii)\,\,\,     a_1=0,  \,b_2=1\,\,  \\[2mm]

\!\!(c) \,\,\,3   \!+\!min\{ a_1, k\!-\! h \}   \!+\!min\{a_2\!-\!a_1 \!-\!1 , k-h\}   \!+\!min\{h\!-\!a_2\!-\!1 , b_3\}   & if & a_1<a_2<h\leq a_3   \\[2mm]
 \hspace{0.5cm}further\,\, \beta_3\geq 4, \,\,  except\, \,the\,  \,cases:\\ 
 \hspace{4cm} (c_1)\,\,   a_1 = a_2 \!-\!1=  min\{h\!-\!2, b_3\}=0   \\[2mm]

\!\! (d) \,\,\,  3 +min\{ a_1, k-h \}+\sum_{i=1,2} min\{a_{i+1}\!-\!a_i\!-\!1 , k-h\}  & if & a_1<a_2<a_3<h  \\[2mm]
\hspace{0.5cm}further\,\, \beta_3=  3\II a_1=0,a_2=1,\, a_3=2  \\[2mm]

  \end{array}\right.$  
  \denu
 \end{lemma}
   Proof. {\it 1.}  Put  $a_1=a,b_1=b$ for simplicity.  
   If $a \geq h$,the induced distinct elements in $C_h$ are: \\$ \{ hn_i ,\,(h-1)n_i+ n_j,\cdots ,hn_j\}$, if $h\leq b$, 
 \\ $ \{ hn_i ,\,(h-1)n_i+ n_j,\cdots ,bn_j\}$, if $h>b$.\\
 Then $\beta_1=1+min\{b , h\}=i+min\{a,b,h,a+b-h\}$.
 \\[1mm]
If $a<h\leq b$, \, then \,  $\{ hn_i  ,\,(h-1)n_i+ n_j,\cdots ,(h-a)n_i+ an_j\}$  are indueced distinct elements of $C_h$. Hence  $\beta_{1}=a+1$ where $a=min\{ a,b,h,a+b-h\}$ since $a\leq a=b-h<b$. 
  \\[1mm]
  If $a,b<h$, \, then the induced distinct elements are: $   an_i +(h-a)n_j,\,(a-1)n_i+(h-a+1)n_j,\cdots ,(h-b)n_i+bn_j\ $, then  $\beta_1=a+b-h+1=k-h+1 $, further $  \beta_{1}= 1+min\{  a,b,h,a+b-h\}=a+b-h(=k-h )$, since $a+b-h<a<h, \,a+b-h<b$.\\[1mm] 
 If $a,b\geq h$, then $\beta_{1 }=1+h$ and   the induced distinct elements are $ \{ hn_i ,\,(h-1)n_i+ n_j,\cdots ,hn_j\}$.  In this case  $| C_h|=h+1$ is maximal. Then, if $a\geq h$, $\beta_1=1+min\{b_1, h\}$.\\[2mm]
 {\it 2}($a$). \, The induced distinct elements $\in C_h$ are $ \,  \{(h-i)n_i+in_j,\,\, 0\leq i\leq min\{b_1,h\}\,\, \}$.\\
{\it 2}($b$).\,   We have $ 1+min\{a_1,k-h\}$ elements $ \in C_h$  induced by $x_1  $, \,  by (1): \\ \centerline{$(a_1-i)n_i+(h-a_1+i)n_j,\,\, \,  0\leq i\leq min\{a_1,k-h\}$}
 Moreover \,  from  $x_2  $    we can extract the following distinct  additional  elements\\ 
 \centerline{$(h-p)n_i+p n_j,\,\, \big(p\geq 0, \, \, h-p\geq a_1+1,\,\, p\leq b_2\big)$ hence with \,$ 0\leq p\leq min\{h-a_1-1,b_2\}$.}\\[2mm]
{\it 2}$(c)$.\, First,  $x_1  $   induces the same elements of $C_h$ considered in ($2.b$);   moreover  from   $x_2 $ one gets other $M$  distinct elements \\
\centerline{$ (a_2-p)n_i+(h-a_2+p)n_j,$ with $ p\geq 0,\,\,a_2-p\geq a_1+1,\,\,h-a_2+p\leq b_2$} (hence $ 0\leq p\leq min \{a_2-a_1-1, k-h\}$),  where 
  \par   either $M=(a_2-a_1)$, \,\, if \, $h-a_1-1\leq b_2$\,\,(elements $(a_2-p)n_i+(h-a_2+p)n_j,$ with $   0\leq p\leq  a_2-a_1-1 $);\par
    or\quad  $M= k-h+1$, \,\,if \, $h-a_1-1 >b_2$
  \, (elements $ \,  a_2n_i+(h-a_2)n_j,\cdots , (h-b_2)n_i+b_2n_j $).
   \\[2mm]\
 {\it 3.} It comes directly by using the same ideas of the proof of  statement ({\it 2}).
 \quad$\diamond$\\[2mm]
This lemma allows to prove the more general 
\begin{prop}\label{supp22} 
Assume  $k\geq 3,\quad 2\leq h\leq k-1$ and $\,\,   C_k \supseteq \{x_1,x_2,\cdots ,x_p\},\, p\leq k+1$, \,\, with\\ 
 \centerline{ $\left\{\begin{array}{llll}x_i = a_i n_i+b_in_j,&1\leq i\leq p, \quad    a_i+b_i=k  \\ 
0\leq a_1<a_2<\cdots <a_p \leq k & (\I  k\geq b_1>b_2>\cdots >b_p\geq 0)  \end{array}\right. $}  and let $\beta_p$ be the number of distinct elements of $C_h$ induced by $\{x_i,\,\,i=1,\cdots ,p\}$. Then $\beta_p\geq min\{ h+1, p\}$, precisely:\\[2mm]
$ \left[\begin{array}{lrllllll}
1+min\{b_1,h\}&if &h\leq a_1<\cdots <a_p\\
    
 i+1 +min\{a_1,k-h\}+  \cdots  +min\{a_i-a_{i-1}-1,k-h\}+min\{h-a_i-1, b_{i+1} \}&if&   <a_i<h\leq a_{i+1}<  &  \\
 && (i\leq p-1)\\
 
 p +min\{a_1,k-h\}+  \cdots  +min\{a_p-a_{p-1}-1,k-h\} &if&  <a_p<h&\\ 

  \end{array}\right.$  
 \end{prop}  

   \begin{lemma} \label{pq}  
    Let   $x_1,x_2$ be distinct elements in $D_{k}$ such that     $  Supp(x_1+e) =\{n_i,n_j\},\quad Supp(x_2+e ) =\{n_t ,n_u  \}$,\,\,$\{n_i,n_j\}\cap\{n_t ,n_u  \}=\emptyset$. \quad  Let$\left\{\begin{array}{l}y_1=x_1+e  = a n_i+b  n_j\\ 
 y_2= x_2+e  = c n_t +d  n_u  \\
 abcd\neq 0,\,\,a+b=k+r_1,\,\, c+d=k+r_2 \end{array}\right.$. \quad \\
Let $r:=min \{r_1,r_2\},$ and let $ 2\leq h\leq k+r $.
 Consider the induced elements  {$z_1=pn_i+(h-p)n_j,\quad z_2=qn_t +(h-q)n_u  \in C_h.$ } 
 \enu \item[] Then   $z_1\neq z_2$ for every $p,q,h$ and $| C_h|\geq  \beta_{ab}+\beta_{cd}$  \, where $\beta_{ab},\beta_{cd}$ are defined in 
 (\ref{supp2}.1) (called $\beta_1$). Consequently\
 \item[] $| C_h| \geq 4 $, \,if $\,\,h<k+r$,\quad
  $| C_h |\geq 3 $, \,if $\,\,h=k+r,\,\, r_1\neq r_2$,\quad
   $| C_h |\geq 2 $,\, if $\,\,h=k+r,\,\, r_1=r_2$.
  \denu
  \end{lemma}
  Proof.  
 Assume $z_1=z_2$, then   by substituting we get\\ $   \left[\begin{array}{llllrr} 
y_1=(a-p)n_i+(b+p-h)n_j+z_1 =(a-p)n_i+(b+p-h)n_j+qn_t +(h-q)n_u  \\
 y_2= (c-q)n_t +(d+q-h)n_u  +z_1=pn_i+(h-p)n_j+(c-q)n_t +(d+q-h)n_u  
\end{array}\right.$.  
 \\
First note that by the assumption, we have $z_1\neq z_2$ if:\,\,$p=q=0,$ or $ p=q-h=0,$ or $ q=h-p=0,$ or $  h=p=q $ ; moreover if three coefficients are $\neq 0$, then $| Supp(y_i)|=3$, against the assumption.    This argument allows to complete  the proof in the  remaining cases which are the following:\\[2mm]
\ $\left[\begin{array}{cccccccccl} 
(a-p)n_i&(b+p-h)n_j&qn_t   &(h-q)n_u  &|&pn_i&(h-p)n_j&(c-q)n_t   &(d+q-h)n_u  &\\
\neq 0 &   0 & 0&\neq 0 &|&\neq 0 &b&c& d-h & \\
 
\neq 0 &   0 & \neq 0&  0 &|&\neq 0 &b&c-h &d & \\
 
  0 &   \neq 0  &  0 &\neq 0 &|&\neq 0 & \neq 0&c& d-h & \\ 
 
    0 &   \neq 0  & \neq 0 &0 &|& \neq 0& \neq 0 &c-h &d&\,   \quad\diamond\\

\end{array}\right.$

  \begin{lemma} \label{3ni}   {  Let $k\geq 3$ and let $x_1,x_2$ be distinct elements in $D_{k}$   such that $| Supp(x_1+e)|=|Supp(x_2+e)|=2$, \par  $\left\{\begin{array}{l} x_1+e  = a n_i+b  n_j\\ 
  x_2+e  = c n_t +d  n_j\\
   \end{array}\right. ,\,  a,b,c,d>0$  with  $\, n_i, n_j,n_t $ distinct elements in $    Ap_1$.  }  \\ For $h<k$, consider the   induced elements in $C_h:$   $   \left\{\begin{array}{llllrr} 
z =pn_i+(h-p)n_j, & 
 z'=qn_t +(h-q)n_j\, 
\end{array}\right\}$. Then   
 \enu
\item $z \neq z'$ in the following cases
  $ \left[\begin{array}{llllrr}
(i)&pq=0  & and& \,\,  p+q>0 \\
 (ii)&
a<h& and&  q\geq max\{1,h-d\}\\
&& or\\
& c<h& and& p\geq max\{1,h-b\}\end{array}\right.$
  \item  In the cases  $\left\{\begin{array}{l} x_1+e  =   n_i+b  n_j\\ 
  x_2+e  =   n_t +d  n_j\\
 
 \end{array}\right.$  we have $| C_h|\geq 3$, \, in the remaining cases  $| C_h|\geq 4.$

   \denu \end{lemma} 
 Proof. {\it 1}. \, $(i)$ is immediate by the assumptions.\\
 $(ii)$. \, It is enough to consider $p>0$ by $(i)$; if $z =z'$, then\\
  $x_2+e=pn_i+(d+q-p)n_j+(c-q)n_t $, where $d+q-p\geq h-p\geq 1$, since $ p\leq a<h$. Then:\\
   $\left[\begin{array}{llllrr}
 |  Supp(x_2+e)|\geq 3& if&c-q>0\\
\quad Supp(x_2+e)= Supp(x_1+e)&if&c-q=0\\ 
 \end{array}\right.$, \quad contradiction in any case.\\[2mm]  
  {\it 2}.\,  We can assume  $0<a\leq c$. The following $z_i$ belong to $C_h$ and are distinct.\\[2mm]
  $\begin{array}{lll}  \left.\begin{array}{lllllllll}
  (i)& a=c=1:  
&z_1 = n_i+(h-1)n_j\\
&& z_2 =h\,\,n_j\\
&&z_ 3=n_t +(h-1)n_j& \vspace{2mm} \end{array}\right.
&\left.\begin{array}{lllllllll} \qquad\qquad  (ii)&  1\leq a, c < h:   & z_1=an_i+(h-a)n_j\\
&( ac\geq 2)  &z_2=(a-1)n_i+(h-a+1)n_j\\ 
& &z_3=cn_t +(h-c)n_j  \\
&&z_4=(c-1)n_t +(h-c+1)n_j&
\vspace{2mm}\\\end{array}\right.\end{array}$\\
$\left.\begin{array}{lllllllll} (iii)&  a<h= c :&z_1= an_i+(h-a)n_j&   \\
&&z_2 =(a-1)n_i+(h-a+1)n_j&\\ 
&&z_3 =c\,\, n_t  & \\   
& &z_4 =(c-1)n_t +n_j&\end{array}\right.
\left.\begin{array}{lllllllll}\vspace{2mm}\\
  (iv)&  a< h < c :&z_1= an_i+(h-a)n_j&\\   
  &&z_2 =(a-1)n_i+(h-a+1)n_j&\\ 
&&z_3 =h\,\, n_t  & \\   
& &z_4 =(h-1)n_t +n_j& \end{array}\right.$\\
$\left.\begin{array}{lllllllll}   (v)  &  h\leq a \leq  c :&z_1= hn_i &\\   
&&z_2 =(h-1)n_i+ n_j&\\ 
&&z_3 =h\,\, n_t  & \\   
& &z_4 =(h-1)n_t +n_j& \vspace{2mm}\\  \end{array}\right.     $  \\
    \noindent The non trivial subcases of $(i),\cdots ,(iv)$ come directly from part {\it 1}.    \\
Case $(v).  $ $z_1=z_4\I $ either $  Supp(x_1+e)=Supp(x_2+e) ($ if $a=h$), or $| Supp (x_1+e)|=3, ($ if $a>h)$, against the assumptions.\\
$z_2=z_3\I $ either   $  Supp(x_2+e)=Supp(x_1+e)\, ($ if $c=h$), or $| Supp (x_2+e)|=3, ($ if $c>h)$, against the assumptions.\quad$\diamond$\\[2mm]
 Next proposition shows that for $k=2$, statement (\ref {pt}.4) holds also at step $k+1$.
   \begin{prop}\label{decr}  
Assume $H_{\!R}$ decreasing at level $\, h$. Then
\enu
\item There exist $x_1,x_2\in D_h$ such that $| Supp(x_r+e)|\geq 2$, \, for \, $r=1,2$.
\item $| C_3 |\geq 4$.
 
\denu
\end{prop}
Proof. {\it 1}.\, Assume $H_{\!R}$ decreasing at level $h$: we know that $  m = | C_h|\leq| D_h|-1$. If $| Supp (x_r+e)|=1 \,\, \forall x_r\in D_h$, then $x_r+e=\alpha_r m_r, \, \alpha_r\geq h+1, m_r\in Ap_1$, with $m_r\neq m_s,$\,\,if $\, r\neq s$. Hence the elements $\{h m_r\} $ would be distinct elements $\in C_h$ and $| C_h|>m$.\\
Hence let $x_1\in D_h$, $| Supp(x_1+e)|\geq 2$, with a maximal representation  \, $x_1+e= \sum  \beta_i m_i $,  with $\beta_i\geq 1, m_i\in Ap_1 $  distinct elements. One induced element $\in C_h $ is \,\,  $c_1=m_i+m_j+ s$ for some $s $ of order $h-2$.  If each $x_i\in D_h, i>1$ has maximal representation of the type  $x_i+e=\alpha_i n_i $,  then $hm_j,\cdots , hm_{m+1}\in C_h$ are distinct elements; further    $hm_i\neq c_1$ for each $i$ (otherwise $ x_i+e= c_1+(\alpha_i-h) m_i$ and so  $| Supp(x_i+e)|\geq 2$). \, Then $|C_h|\geq m+1$, contradiction.\\[2mm]
{\it 2}.\,  By (\ref{pt}.3 and 4)  we know that  if $h>2$ then $| C_j|\geq j+1$ for all $2\leq j\leq h$. It remains to prove the statement for $h=2$, \,\, and \,\, $| D_2|\geq 4$ by \cite[prop.2.4]{ddm}. \\
Since $C_3=Ap_3\cup  (D_2^3+e) $, the fact is true if $ord(x+e)=3\,\, \forall x\in D_2$.\,\, Hence we assume there exists $x\in D_2$ such that $ord(x+e)\geq 4$.\\
In this proof, for simplicity, we denote the elements  of $Supp(x)$ with $m_1,m_2,m_3\cdots $ and assume $m_1<m_2<m_3$. 
 \\
If  $| Supp(x+e)|\geq 4: x+e=m_1+m_2+m_3 +m_4 +s$, with $ m_1<m_2<m_3 <m_4 $, then  $| C_3|\geq 4$  by (\ref{suppq}.$a$) .
\\[2mm]
If $| Supp(x+e)|=3, \,x+e=\beta_1m_1+\beta_2m_2+\beta_3m_3$   \,with   $\beta_i\geq 1,\,  \sum  \beta_i\geq 5$),we find  $4 $ induced elements in $\ C_3$ as follows: \par
when $\beta_1\geq 2, \beta_2\geq 2$: $\{m_1+m_2+m_3, 2m_1+m_2, 2m_1+m_3,2m_2+m_3\}; $\par
when $\beta_1\geq 2, \beta_3\geq 2$: $\{m_1+m_2+m_3, 2m_1+m_2, 2m_1+m_3, m_2+2m_3\}; $\par
when $\beta_1, \beta_2\geq 2,\beta_3\geq 2$: $\{m_1+m_2+m_3, 2m_2+m_3, m_1+2m_2, m_2+2m_3\}; $\par
when  $\beta_1\geq 3$: $\{m_1+m_2+m_3, 2m_1+m_2, 2m_1+m_3, 3m_1\}. $
\\[2mm]
If $| Supp(x+e)|=3,  \,x+e=\beta_1m_1+\beta_2m_2+\beta_3m_3$ with $\beta_i\geq 1,\,  \sum  \beta_i=4$, then  $x+e$  induces in  $C_3$    three distinct elements $ c_1,c_2,c_3$ (see the following table). To get the fourth element, we take $y\in D_2, y\neq x,  | Supp(y+e)|\geq 2$ (by 1). If  $Supp(x+e)\neq Supp(y+e) $  and each element induced  by $y+e$ in $C_3$ is also induced by $x+e$, then  \, we can always reduce  to the case $Supp(x+e)\supseteq Supp(y+e)$, by suitable substitutions (in one or more steps).\\  Hence we can assume $Supp(x+e)\supseteq Supp(y+e)$; we study in the next table one of the possible subcases, the remaining are similar. Let \,\\[1mm]\centerline {$\left\{\begin{array}{llllllcc}
    x+e= m_1+ 2m_2+ m_3\\
y+e=\beta'_1m_1+\beta'_2m_2+\beta'_3m_3,    &\quad \sum \beta'_i \geq 3 
 \end{array}\right.$ } \\[2mm]   
    If  $\sum \beta'_i = 3$, then $y+e\neq z$ for each $z\in C_3$ induced by $x+e$, otherwise $ord(x)>ord (y)$. Hence assume  $ \sum \beta'_i > 3$;
   \par then  $ C_3\supseteq\left[\begin{array}{llllllcc}
c_1 =m_1+2m_2\\
c_2 =m_1+m_2+m_3 \\
c_3 =2m_2+m_3 \\ 
 and,  \,\,according \,\,to\,\,the\,\,values\,\,of \, the \,\,\beta'_i,  \\
c_ 4=2m_1+m_2&if& \beta_1'=2,\beta_2'\geq 1 \\

c_5=m_2+2m_3   &if& \beta_1'\leq 1,1\leq \beta_2'\leq 2,\beta_3'= 2 &\\
c_ 6=3m_3 &if& \beta_3'\geq 3 \\

c_7 =3m_1&if& \beta_1'\geq 3 \\
c_8 =3m_2, \, or\, \, 2m_1+m_2, \, or\, \, m_2+2m_3 &if& \beta_2'\geq 3\\
 c_ 9=2m_1+m_3   &if& \beta_1'=2,\beta_2'=0,\beta_3'\geq 1 \\
\end{array}\right.$\\ [2mm]
Clearly the elements $c_i,\,i=4,\cdots ,7$ are distinct from $ c_1,c_2,c_3$. 
In  case $\beta_2'\geq 3$ we have either $y+e=3m_2+m_i$, \,$(i=1$ or $3$)  \,\, or \,\,$y+e=3m_2$ . In the first case if $\, 3m_2=m_1+m_2+m_3 $, then
$2m_2= m_1+m_3  \I y+e=  2m_1+ m_2+ m_3 , $\, or \, $y+e= m_1+ m_2+ 2m_3 $ and   we can add $c_4$ \, resp. \, $c_5$.\\
In  case $y+e=3m_2 $, we get  $ 3m_2\notin\{c_1,c_2,c_3\}$, otherwise     
$2m_2= m_1+m_3  \I y+e=   m_1+ m_2+ m_3 , $ i.e., $x=y+m_2$, impossible since $ord(x)=ord(y)$.\\
If  $\beta_1'=2,\beta_2'=0,\beta_3'\geq 1$, then $c_9\notin\{c_1,c_2,c_3\}$ otherwise $2m_1+m_3 =m_1+2m_2\I m_1+m_3 = 2m_2\I x-y=(2-\alpha)m_3 ,$ impossible in any case.\\[2mm]
   Finally we assume:  $| Supp(x+e)|=| Supp(y+e)|=2,\,\, ord(x+e)\geq 4$.
By lemmas (\ref{pq}) and (\ref{3ni}), it remains to analyse the cases\\ 
\centerline{$(I)\left\{\begin{array}{l} x +e  =   m_1+b  m_2\\ 
 y+e  =   m_3 +d m_2\\
  b\geq 3, \,d\geq 2 
  \end{array}\right.$ \qquad  $(I\!I)\left\{\begin{array}{l} x_1+e  =  a  m_1+b  m_2\\ 
  x_2+e  =  c m_1+d m_2\\
 a+b\geq 4, \,c+d\geq 3  
  \end{array}\right.$  \qquad $(I\!I\!I) \left\{\begin{array}{lll} x +e  =  a m_1+b  m_2&a+b=3\\ 
 y+e  =   cm_3 +dm_4 &c+d\geq 4\\
  \end{array}\right.$}
 $(I)$. \, The induced distinct elements are $\left\{\begin{array}{llllllcc}
c_1 =m_1+2m_2&
c_2 =3m_2&
c_3 =2m_2+m_3 \\ 
 \end{array}\right\}$. \\
  By assumption,   there exist two other elements $z_1,z_2\in D_2$;   by the above tools and by (\ref{short}),  we can restrict to the cases \\
  \centerline{   $ b\geq 3, \,d\geq 2 :\left[\begin{array}{llll} 
(i)&z_1+e  =  \alpha m_1+\beta m_2&\alpha\neq 1&\alpha  +\beta \geq 3 \\ 
 (ii)&z_1+e  =   \alpha m_3 +\beta m_2&\alpha\neq 1&\alpha  +\beta \geq 3\\
 (iii)& z_1+e  =   \alpha m_ 1+\beta m_3 &\alpha  +\beta \geq 3     \end{array}\right.$.}   
  In the first two cases, if $ \alpha >0,$ then $\alpha\geq 2$ (by \ref{short}.2), therefore if $\beta>0$  we obtain $| C_3|\geq 4$ by  applying (\ref{3ni}.{\it 2}) to the pair of elements $\{z_1+e,\, y+e\}$ (resp $ \{x+e,\,z_1+e \}$). In case $(iii)$, again,  from (\ref{3ni}.{\it 2}),     by   substituting in $\{x+e,\,z_1+e \}$ $m_1'=m_2, m_2'=m_1, m_3 '=m_3 $, we deduce $| C_3|\geq 4$, \,  if $\,\alpha\beta>0$.
  The remaining possibilities are\\ \centerline{ $\left.\begin{array}{llll}
 (i')&z_1+e  =   \beta m_2&  \beta\geq 3 \\ 
 (ii')&z_1+e  =   \alpha  m_1&  \alpha \geq 3  \\ 
 (iii')&z_1+e  =    \alpha m_3  & \alpha\geq 4& 
  \end{array}\right. .$}\\ In case $(i')$ we can assume  the element $z_2\in D_2$,
  verifies $z_2+e= \alpha m_1,\,\, \alpha\geq 3$ \,\,  (or $z_2+e= \alpha m_3 $). \\
  If $3m_1=2m_2+m_3 $, then, either $| Supp(z_2+e)|=3$ and we are done, or 
  $z_2+e =2m_2+m_3 $, impossible   since $y+e=m_3 +dm_2$. \\
  Similarly we can   solve case $(ii')$.\\
   In case $(iii')$, if $3m_3 =m_1+2m_2$, \, then   \, $z_1+e=m_1+2m_2+(\alpha-3)m_3 . $ The case $\alpha\geq 5$ has already been proved above. If $\alpha=4,$ we consider $m_1+m_2+m_3 \in C_3$ which, in this situation,  is distinct from $c_1,c_2,c_3$.\\[2mm]
  ($I\!I)$. If we have $4$ elements with $support\subseteq\{m_1,m_2\}$, then we are done by (\ref{supp22}). The other cases can be studied among the ones of shape $(I)$ or $(I\!I\!I)$.  \\[2mm]
  $(I\!I\!I)$. \,  By \ref{pq}, it remains   to consider elements $z_i+e=\alpha m_1+\beta m_2$,\, with $\alpha+\beta=3$. In this case,   $C_3\supseteq\{x+e,z_1+e, \}\cup\{y_1,y_2, \,$ induced by $y+e\} $.
  \quad$\diamond$
\section{Structure of  \, $  C_3$.}

 \begin{thm}\label{c23}  
 Assume \,  $| Ap_2|=3$.  \, Then \enu

 \item If $| C_3|\geq 2$,\, then $| Supp(x)|\leq 2$,  \,for each $x\in C_3$.
\item  If $| C_3|\geq 4$, \, then 
\begin{enumerate}  
\item   There exist $x_1,x_2\in C_3$ such that $| Supp(x_i)|=2$.\\[-5mm]
\item   There exist $n_i,n_j\in Ap_1$ such that $ C_2(=Ap_2)=  \left\{\begin{array}{lllll}
 2n_i \\
 n_i+n_j,\\
2n_j\\
 \end{array}\right.\quad C_3= \left\{\begin{array}{lllll}
  3n_i \\
  2n_i+n_j\\
  n_i+2n_j\\
  3n_j
 \end{array}\right.$.
  \denu
\denu
\end{thm}
Proof. {\it 1}.\,  If $x\in C_3$, then $| Supp (x  )|\leq 3$. First we show that if $x_1,x_2\in C_3$, then $| Supp (x_i  )|\leq 2$ for $i=1,2$.
Assume that $x_1  = n_i+n_j+n_k,\,\, n_i<n_j<n_k$.\,  Then, by (\ref{PT0}.1$a$)  and the assumption, we deduce that\\
\centerline{$Ap_2=\{ n_i+n_j ,n_i+ n_k, n_j+n_k\}\qquad (*)$}
Hence $2n_i, 2n_k\notin Ap$,\,\, $2n_j\in Ap\II n_i+n_k=2n_j$.
If $x_2\in C_3 , x_2\neq x_1$ there are three cases:\enua 
\item If \, $| Supp(x_2  )|=1,$ hence  $x_2  =3n_t$, then $2n_t\in C_2$, by (\ref {PT0}.1$a$). Hence by $(*)$ $2n_t=\left[\begin{array}{ccc}
n_i+n_j\\
 n_i+n_k\\
 n_j+n_k
 \end{array}\right.$\\ In every case we see that $| Supp(x_2  )|\geq 2$ ,  a contradiction.
 \item If \, $| Supp(x_2  )|=2,$  $x_2  =2n_t+n_v, t\neq v$, then $2n_t\in C_2$, by (\ref {PT0}.1$a$).\\
  Hence by $(*)$ $2n_t=\left[\begin{array}{lllll}
n_i+n_j&\I x_2  =n_i+n_j+n_v \I n_v\in\{n_i,n_j\}\\
&n_v=n_i\I x_2  =2n_i+n_j  \I 2n_i\in Ap_2\\
&n_v=n_j\I x_2  = n_i+2n_j \I n_i+n_k=2n_j\in Ap_2 \I 2n_i\in Ap_2\\
 n_i+n_k &({\rm similar})\\
 n_j+n_k& ({\rm similar})
 \end{array}\right.$ \\ every case contradicts equality $(*)$.
 
\item If $\, | Supp(x_2  )|=3,$ \, $x_2  = n_t+n_v+n_w, \,\, n_t<n_v<n_w  $, \, by \, $(*)$ we must have: $n_t+n_v=n_i+n_j$, hence $x_2  =n_i+n_j+n_w$. This equality implies that \\
$\left\{\begin{array}{lllll}
n_w\neq n_k &({\rm since}\, x_2\neq x_1)\\
n_w\neq n_j&({\rm otherwise }\,\,\, x_2  =n_i+2n_j=2n_i+n_k,\quad {\rm  see}(*))  \\
n_i+n_w\in Ap_2&\I n_i+n_w=n_j+n_k
 \end{array}\right.$. \\[2mm]
 Hence $x_2  =2n_j+n_k=n_i+2n_k$, contradiction. \denu
{\it 2}($  a$). Let $x_i\in C_3 , 1\leq i\leq 4,$  distinct elements such that  $1\leq | Supp (x_i  )|\leq 2$ (by {\it 1}). If $| Supp (x_i  )|=1$, for $i=1,2,3$, $\left\{\begin{array}{lllll}
x_1  =&3n_a\\
x_2  =&3n_i\\
x_3  =&3n_c\\
 \end{array}\right.$,  then by (\ref {PT0}.1$a$) $C_2=Ap_2=\{2n_a,2n_i,2n_c\}$. Since $| Ap_2|=3$, by ({\it 1}) and  (\ref{PT0}.1$a$), we have   $| Supp (x_4  )|=2$:   we can assume $x_4  = 2n_a+n_d, n_a+n_d=2n_i$. Hence $ x_2  =n_a+n_d+n_i$, with $n_a\neq n_d$ and so $| Supp (x_2  )|\geq 2$, against the assumption.\\[2mm]
  {\it 2}$(b)$. Let $\{x_1,x_2,x_3,x_4\}\subseteq C_3$; by {\it 2}({\it a}),   we can assume that   $\left\{\begin{array}{lllll}
x_1   = 2n_i+n_j,& n_i\neq n_j\\

x_2   =2n_h+n_k& n_h\neq n_k\\
 \end{array}\right.$. \\ We   prove that $| Ap_2|=3\I  n_i=n_k,n_j=n_h$, \, i.e. $Supp(x_1   )=Supp(x_2   ), x_2=n_i+2n_j$.\\
Note that  $| Ap_2|= 3\I |\{ n_i,n_j,n_h,n_k\}|\leq 3$, otherwise, $\{ 2n_i, 2n_h,n_i+n_j, n_h+n_k\} \subseteq Ap_2\I | Ap_2|\geq 4$. In fact if $2n_i=n_h+n_k\I | Supp(x_1+e)|\geq 3$ ( the other cases are similar). Hence  there are three possible distinct situations: \\[2mm] 
 \centerline{ ({$ b_1$}) $  \left\{\begin{array}{lllll}
x_1   = 2n_i+n_j,& n_i\neq n_j, \\ &n_i\neq n_h\\
x_2   =2n_h+n_j& n_h\neq n_j\\
\end{array}
\right.,$     \quad  ({$ b_2$})   $\left\{\begin{array}{lllll}
x_1   =  2n_i+n_j,& n_j\neq n_i\\
x_2   =2n_h+ n_i& n_h\neq n_i\\
\end{array}
\right.,  $ \quad  ({$ b_3$})   $\left\{\begin{array}{lllll}
x_1   =  2n_i+ n_j,& n_i\neq n_j\\
x_2   = 2n_i+n_k& n_k\neq n_j\\
\end{array}
\right..$ }\\[2mm] 
  ({$ b_1$}).\,   We have $Ap_2\supseteq \{2n_i,n_i+n_j, 2n_h,n_h+n_j\},$ hence $| Ap_2|=3\I  $ \\ \centerline{ $\left[\begin{array}{lllll}
 either&n_i+n_j=2n_h     \I x_2   =n_i+2n_j& ({\rm thesis})\\
 \quad or& 
 2n_i=  n_h+ n_j&({\rm analogous})\,.\\
\end{array}\right.$  }\\[2mm]
({$ b_2$}).\,   We have  \,\,$Ap_2\supseteq \{2n_i,n_i+n_j, n_h+n_i, 2n_h\} \I 2n_h=n_i+n_j\I $ either $ x_2=x_1$\\[2mm] (against the assumption), or \, $n_h=n_j (\I x_2   =2n_j+n_i)$.\\[2mm]
 ({$ b_3$}).\,  This case cannot happen. In fact we have  $Ap_2= \{2n_i,n_i+n_j, n_i+n_k\}  $. By similar arguments as above one can see that for any other element  $x\in C_3,\,  x\neq x_1,x_2$, the  maximal representation of $x $ must be written   as $x=3n_i $. In fact the other possible representations are incompatible; for instance,   $x=3n_j , $  with $2n_j=n_k+ n_i\I x=n_i+n_j+n_k$, impossible by (2). This would mean that $| C_3|\leq 3$, against the assumption.\\ According to the above facts we deduce that $C_3=\{3n_i ,
  2n_i+n_j,
  n_i+2n_j,  3n_j\}$. 
  \quad$\diamond$ 
   \begin{ex}\label{es2}  {\rm   According to the notation of (\ref{c23}), we show several examples of semigroups which verify the assumptions of (\ref{c23}):    $ Ap_2=\{2n_i,n_i+n_j,2n_j\}$. The first example shows that the conditions (\ref{c23}.2$b$) are, in general, not sufficient to have $H_R$ decreasing.}
  \enu 
    \item  {\rm Let $S=< 19,21,24,  47,49,50,51,52,53,54,55,56,58,60>$ and let $n_i=21,n_j=24$. Then:   $Ap_2=\{ 42,45,48\}, $ \, $ Ap_3=\{63\},Ap_4=\{84\}$,  $\,v=e-5$,  \,\,  $H_R=[1,14,14,14,16,18,19\rightarrow]$ is non-decreasing. }
 
  \item {\rm Let $S=< 19,21,24,65,68,70,71,73,74,75,77,79>$ with $n_i=21,n_j=24$. 
 Then: \\$\ell  =2 ,   Ap_2=\{42,45,48\},\,\, Ap_3=\{3n_i,2n_i+n_j,n_i+2n_j,3n_j\} ,\,\,v=e-7$,\, \, $D_2+e=\{ 4n_i, 3n_i+n_j,...,n_i+3n_j\},\, D_2=\{65,68, 71,74, 77\}$\quad $H_R$ decreases at level 2,\, $H_R=[1,{\bf12,10},11,15,18,19\rightarrow]$.}
  
   \item {\rm Let $S=< 19, 21, 24, 46, 49, 51, 52, 54, 55, 56, 58, 60>$ with $n_i=21,n_j=24$. Then:\\  $\ell  =3,  \, | Ap_3|=4 ,\,\,v=e-7$,\, $H_R$ decreasing at level 3,\,$H_R=[1,12,{\bf15,14},16,18,19\rightarrow]$. }
   \item    {\rm $S=< 30,33,37,64,68,71,73,75, 76,78,79,80,82\rightarrow 89,91,92,95>$ with $n_i=33,n_j=37$. Then:\\ $   Ap_2=\{2n_i,n_i+n_j,2n_j\},   \, Ap_3=\{3n_i,n_i+2n_j,3n_j\},\, Ap_4=\{132=4n_i\},  \,\,v= e-7 $,\,  $H_R$ decreases at level 4: $H_R=[1,23,25,{\bf25,24},27,28,29,30\rightarrow]$.}
        
 \denu 
 \end{ex}
     \begin{lemma} \label{preapj1}  Let $d =max\{ord(\sigma) \,|\,\, \sigma\in Ap\}$. Assume     there exists $3\leq r\leq d $   such that $\, | Ap_r|=1.$   Let $\, r_0:=min\{j\,\,|\,\, | Ap_j|=1\}  \,\,(r_0\geq 3)$,\quad then$:$
 \begin{enumerate} 
\item   If $r_0<d $, there exists    $n_i\in Ap_1$ such that $Ap_k=
\{kn_i\}$  \,\, for \, $r_0\leq k\leq d  \, ($ and \, $kn_i\in Ap_k$, \, for $k<r_0)$. 
\item If there exists   $ g\in D_k,\,( k\geq 2)$ such that $ord(g+e)=k+p,\, $ with $p\geq r_0-1$, then  
  \, $Ap_d = d  n_i \, (n_i\in Ap_1)$,  \,and\,  $g+e=(k+p)n_i>d n_i ;$  \ further  such element  $\,g\,$ is  unique. 
  
\denu   
 \end{lemma}
 Proof. {\it 1}. By the Admissibility Theorem of Macaulay for     $H_{R'}$, we know that $| Ap_k|=1$, for $r_0\leq k\leq d $. Now suppose $r_0<d $, and $| Supp(\sigma)|\geq 2, \,\,\sigma\in Ap_d $. If  $n_i,n_j\in Supp(\sigma)$, then, by (\ref{PT0}), the elements $\sigma-n_i,\,\sigma-n_j$ would be distinct elements in $Ap_{d -1}$,\, contradiction.\\[2mm]
 {\it 2}. By the assumption and (\ref{CnonAp}.1$\,c$), we have $r_0\leq p+1\leq d $ \, and $\, | Supp(g+e)|\leq | Ap_{p+1}|=1$.  Therefore $g+e=(k+p)n_i$, with $n_i\in Ap_1$. Then for $r_0\leq j\leq p+1$, the induced element   $jn_i\in C_j$ belongs to $Ap_j$, by (\ref{CnonAp}.1\,$a$). Then $Ap_j=\{jn_i\}$ for $r_0\leq j\leq d $. Hence $k+p>d $. For each $k$, if such $g$ exists, then it is unique by (\ref{short}.2).  
 \quad$\diamond$

  \begin{prop} \label{apj1}   Assume $| Ap_2|=3$, $| Ap_3|=1$  and    $H_R$ decreasing.  Let $\ell =min\{h\,| \, H_R$ decreases at level  $h\}$ and let $d =max\{ord(\sigma) \,|\,\, \sigma\in Ap\} $ . We have: \enu
  \item $\,\,  \ell \leq d $,\,   there exist \,  $ n_i,n_j \in Ap_1$,   such that   $(d +1)n_i\in D_{\ell }+e$,   \\  $C_2 =  \left\{\begin{array}{lllll}
 2n_i \\
 n_i+n_j,\\
2n_j\\
 \end{array}\right.  \,C_3= \left\{\begin{array}{lllll}
  3n_i  \\
  2n_i+n_j\\
  n_i+2n_j\\
  3n_j
 \end{array}\right., \cdots,\, C_{\ell }= \left\{\begin{array}{lllll}
  {\ell }n_i \\

  ({\ell }-1)n_i+ n_j\\
   ({\ell }-2)n_i+2 n_j\,,\\
  \cdots \\
  {\ell } n_j
 \end{array}\right. $\\ 
 if \, $r<\ell$, \quad  $C_{r+1}=(D_r+e)\cup Ap_{r+1}$,\\[2mm]
 $  D_{\ell }+e= \left\{\begin{array}{lllll}
 ( d +1) n_i ,\,
\ell n_i+n_j  ,\,
 (\ell -1) n_i+2n_j,  
 \cdots  ,\,
  (\ell +1)n_j 
\end{array}\right\}$   

and,   \, if \,  $(\ell ,d )\neq (3,3)$,\, then   \, $Ap_k=\{kn_i\}$, \,  for all $ \,\, k\in[3,  d ] $.
 \item  The semigroup  $  S$ is not   symmetric.
    \denu
 \end{prop}
Proof. {\it 1}.\,   By (\ref{decr}.2) and (\ref{c23}.2\,$b$) we get the structure of $C_2,C_3$; also,  $| C_r|\leq r+1$ for each $r\geq 0$, since $Supp(C_r)\subseteq \{n_i,n_j\}$ by (\ref{PT0}.1$b$). Further  if $H_R$ decreases at level 2, then $| (D_2+e)|\geq 4$, hence there exists $g\in D_2$ such that $ord(g+e)\geq 4$  (since $|  C_3\setminus Ap_3 |=3)$: by (\ref{preapj1}.2) with $r_0=3,\,p\geq 2$,  we get   $g+e=\alpha n_i$. \\
If $H_R$ decreases at level  $ {\ell }\geq 3$, then  for any $3\leq r\leq {\ell } $, \,\, $| C_r|\geq r+1,$ by  (\ref{pt}.4).
Hence, for $r\leq {\ell } $:  
\\ 
\centerline{\ $C_r=\{rn_i,(r-1)n_i+n_j,\cdots , rn_j\},\,\, | C_r|=r+1.$}
 If $r<\ell\, (\leq d )$, then $| D_r|\leq | C_r|=r+1$, and $r+1=|C_{r+1}\setminus Ap|=  |D_r^{r+1}+e | $; hence $D_r^{r+1}=D_r$.  Further  \\ \centerline{ \quad  $C_{\ell +1}=\{Ap_{\ell +1}\}\cup \big( D_{\ell  }^ {(\ell +1)}+e\big) $.\qquad($*$)}  
 Now, when $ \ell < d $, there exists $g\in D_{\ell }$ such that $ord(g+e)> {\ell }+1$. Hence $g+e=\lambda n_i$ and $Ap_d =d n_i$ by (\ref{preapj1}.2).\quad
If ${\ell }=d $, then \,\,  either $  D_d +e\subseteq C_{ d +1}$ and $| C_{d +1}|= d +2$\I  
$(d +1)n_i=g+e$ with $g\in D_d $, \\
or \,\,there exists $g\in D_d , ord(g+e)>d +1$ and $g+e=\alpha n_i$, hence \,$Ap_d =d n_i$ by (\ref{preapj1}.2) as above. \\
Now we show that ${\ell }\leq d $.\quad If ${\ell }\geq d +1$, then    $(d +1)n_i\in C_{d +1}$; but $(d +1)n_i\notin Ap\I (d +1)n_i=d+e$, with $d\in D_k,\, k\leq d  $, hence $H_R$ decreases at level  $ \leq  d .$\\ 
Finally we prove that $(d +1)n_i\in D_{\ell }+e$: we already know that there exists  $d\in D_{\ell }$ with  $d+e=\alpha n_i\notin Ap$. Then $\alpha \geq d +1$, since $\alpha n_i\in M+e$. If   $\alpha > d +1 $, then  $ord((d +1)n_i-e)< ord (\alpha n_i -e)=\ell -1$, i.e. $ (d +1)n_i-e\in D_k$, with $k< \ell $, impossible by ($*$). \\[2mm] 
{\it 2}.\, Assume $S$ symmetric: then it is well known that for each $n_{\alpha}\in Ap$ there exist $n_{\beta}\in Ap$ with $n_{\alpha}+n_{\beta}= e+f   $ ($f$ is the Frobenius number). Clearly, $ e+f \in Ap_d $. If $e+f  =d n_i$,    in particular there exist $n_r,n_s\in Ap$ such that $n_j+n_{r}=(n_i+n_j)+ n_{s}=d n_i$   (because $C_2\subseteq Ap $).\quad Hence \\
$    n_i+n_s=n_r\in Ap\I \left[\begin{array}{lllll}
  either &n_s=n_j\I   2n_j+n_i=d n_i\notin D_2+e,\,\,   impossible\\
  or&n_s=\lambda n_i (\lambda<d -1)\I n_j=(d -\lambda-1)n_i,\,\,   impossible\end{array}\right.$.\\ If $(\ell ,d )=(3,3)$, we can have $e+f=2n_i+n_j$; in this case, since $3n_i+n_j-\mu e\in Ap_1$ ( with $\mu\in\{1,2\}$ ), there exists
  $n_r$ such that $3n_i+n_j-\mu e +n_r=2n_i+n_j $, hence $n_i+n_r=\mu e$, impossible.   $\quad\diamond$ \\
 
\begin{ex}\label{esj1}  {\rm   According to the notation of (\ref{apj1}), we show several examples of semigroups with $v=e-4$, or $v=e-5$, which verify the assumptions of (\ref{apj1}):    $   Ap_2=\{2n_i,n_i+n_j,2n_j\}$, \, $|  Ap_3|=1$\,\, (see Proposition \ref{apj1}.1   and next  Theorem \ref{31}). In particular examples 2 and 3 show that, in case $ \e=d  =3$, we can have both $Ap_3 = \{3n_i\}  $ and $Ap_3\neq \{3n_i\} $}
 
  \enu
  \item {\rm Let $S=<17,19,22,43,45,46,47,48,49,50,52,54,59>$ and let $n_i=19, n_j=22$.  
Then $\,\,v=e-4$,\\
  $Ap_2= \{38,41,44\}$, $Ap_3=\{3n_i\}= \{57\}$,  \,\,$D_2=\{43=2n_i+n_j-e,\,46=n_i+2n_j-e ,\,49=3 n_j-e, \, 59=4n_i -e\}$; $\e=2, \, d=3,$\,\, $H_R= [ 1,{\bf 13, 12}, 13, 15, 16, 17 ]$.}
 
  \item {\rm Let $S=< 19,21,24, 46,47,49,50,51,52,53,54,55,56,58,60>$ and let $n_i=21,n_j=24$. Then $\,\,v=e-4$,\\
  $Ap_2=\{ 42,45,48\}, Ap_3=\{63(=3n_i)\},$.  
     $D_2+e= \{2n_i+n_j=66, n_i+2n_j=69,3n_j=72\}$, $C_3=(D_2+e)\cup\{63\}$, $\,D_3+e= \{4n_i,3n_i+n_j,2n_i+2n_j, n_i+3n_j,4n_j\}$;\,\, $\ell  =d =3, $\,\, $H_R=[1,15,{\bf15,14},16,18,19\rightarrow]$. }
     
 \item {\rm Let $S=< 19,21,24,44,46,49,50,51,52,53,54,55,56,58,60>$ and let $n_i=21,n_j=24$.  Then $\,\,v=e-4$,\\
  $Ap_2=\{ 42,45,48\}, Ap_3=\{66(=2n_i+n_j)\} \neq \{3n_i\}$, \,
   $D_2+e=\{3n_i=63, n_i+2n_j=69,3n_j=72\}$,  $C_3=(D_2+e)\cup\{66\}$; \quad $\ell  =d =3 $,   \quad $H_R=[1,15,{\bf15,14},16,18,19\rightarrow]$. }
   
  \item {\rm $S=< 30,33,37,73,75\rightarrow 89, 91,92,94,95,98,101>$ and let $n_i=33,n_j=37$. Then: \\
  $Ap_2=\{ 66,70,74\},$  $ Ap_3=\{99\}=\{3n_i\},\,Ap_4=\{132\}=\{4n_i\}$, \, 
   $ v=e-5$. \, One can check that the subsets $C_i,D_i$, for $i\leq 4$ have the structure described in (\ref{apj1}); \,\,   $ \ell  =d =4$,\,\,$H_R=[1,25,25,{\bf25,24},25,27,29,30\rightarrow]$. }
 \denu 
 \end{ex}
 
 {\begin{thm} \label{ap24} 
 Assume the Hilbert function $H_{\!R}$ decreasing     and  \, $| Ap_2|=4$. Then one of the following cases holds       $(n_i,n_j,n_k,n_h$ denote distinct elements in $ Ap_1):$ 
 \enua 
 \item $  Ap_2=\{  2n_i, n_i +n_j ,n_i +n_k , n_j+n_k   \} , \quad   C_3=\{n_i+n_j+n_k, 3n_i, 2n_i+ n_j,  2n_i+n_k\}$. 
\item   $Ap_2= \{ 2n_i ,n_i+n_j, 2n_j, n_i+n_k \},\quad C_3\subseteq\{3n_i,2n_i+ n_j,n_i+2n_j, 3n_j, 2n_i+n_k\}$.
 \item   $Ap_2= \{2n_i,  n_i+n_j ,2n_j,n_h+n_k\}, \quad C_3=\{3n_i, 2n_i+n_j,n_i+2n_j,3n_j \} .$  
 \item$Ap_2= \{2n_i,n_i+n_j,2n_j, 2n_k \}, \quad  C_3\subseteq \{3n_i,2n_i+ n_j,n_i+2n_j,3n_j,3n_k\}$.  
 
\item $Ap_2= \{2n_i ,2n_j, n_i+n_k,n_j+n_k \},\quad C_3=\{3n_i,2n_i+n_k,2n_j+n_k,3n_j\}$.

\denu 

 \end{thm}}
 Proof. \,  Clearly  $| Supp(Ap_2)|\geq 3$, since $| Ap_2|=4$. Assume $H_R$ decreasing. By (\ref{decr}.1-2) and by (\ref{PT0}.1\,$a$),  $| C_3|\geq 4$ and there exist at least two elements $x_1,x_2\in C_3$ with $| Supp( x_i)|\geq 2$.  We proceed step-by-step.  First we assume that    $ | Supp(x_1)|=3,\, \,x_1=n_i+n_j+n_k$.\quad  \\[2mm] 
 Step {\it 1.}\, Let   \, $x_2= n_p+n_q+n_r$.
 Then we get the following induced elements   $\in Ap_2:\{ n_i +n_j ,n_i +n_k ,n_j +n_k ,n_p +n_q, n_p+n_r,n_q +n_r \}$\,    (\ref{PT0}.1\,$a$); \,
  $| Ap_2|=4\I n_p +n_q=n_i+n_j\I
  x_2=n_i +n_j+n_r$ therefore the induced elements in $Ap_2$ can be written as $ \{ n_i +n_j ,n_i +n_k ,n_j +n_k ,  n_i+n_r,n_j +n_r \}$. Again we deduce that  either $n_i+n_r=n_j+n_k$, i.e. $x_2=2n_j+n_k$, \,\,  or $n_r\in\{n_i,n_j\}$. In any case \,    there exists a maximal representation of $x_2$, with $Supp( x_2)\subseteq Supp( x_1)$. Then we can assume that there exists at most one element $x $ in $C_3$ with $ |  Supp(x)|=3$.\\[2mm] 
 Step {\it 2.}\, Let $  \,x_2=2n_p+n_q $. Then we get the following induced elements   $\in Ap_2:\{ n_i +n_j ,n_i +n_k ,n_j +n_k ,n_p +n_q,2n_p\}$, $| Ap_2|=4\I   n_i +n_j=n_p +n_q,$ or \,\, $n_i +n_j=2n_p  $, and so, by using   Step {\it 1},\, we obtain $Supp( x_2)\subseteq Supp( x_1)$.  \\[2mm]  
 Step {\it 3.}\,   Now assume $| Supp(x)|=
\leq 2$ for each $x\in C_3$, then  $| Supp(x_1 )|=  | Supp(x_2 )|= 2. $ \quad If $Supp(x_1 )\cap Supp(x_2)=\emptyset, \quad      x_1=2n_i+ n_j ,\quad 
  x_2= 2n_h+ n_k   ,$ \,then \\ \centerline{ $Ap_2= \{2n_i,n_i+n_j, 2n_h,n_h+n_k\}$,} where these elements are distinct, according to the assumptions. Then   $C_3\setminus\{x_1,x_2\} \subseteq \{    3n_i,3n_h\}$. In fact, by ({\it 1}), any other possible choice   $x \in C_3$ with $| Supp(x )|\leq 2$,  $x$ distinct from 
$3n_i,3n_h$ contradicts some of the assumptions, for example :\\
$x =2n_i+n_k \I n_i+n_k=2n_h\in Ap_2$, impossible since $\I x_2=n_i+2n_k\I n_i\in Supp(x_1 )\cap Supp(x_2)$,\\
 $x = 3n_p, \,p\neq i,h\I 2n_p=n_i+n_j$, impossible since $\I| Supp(x )|=3$.  Hence \\ \centerline{$ C_3=\{3n_i, 2n_i+n_j, 3n_h,2n_h+n_k\},\quad | C_3|=4$.}
Note that the conditions $H_R$ decreasing at any level $k$, $| Ap_2|=4,\ | \geq 4$,  imply $  |  D_k|\geq 5$ (by the assumptions and, for $k\geq 3$, by (\ref{pt}.4)). Since for every $y \in D_k+e$ we have $y \in C_h,\, h\geq k+1\geq 3$  and so $y $ induces elements $\in C_3$, we deduce that these conditions are incompatible. Hence $| Supp(C_3)|\leq 3$.\\[2mm]  
 Step {\it 4.}\, Now we consider the   situation  $| Supp(x_1)|=3$.\\
   Let $ n_i,n_j,n_k$ be distinct elements in $Ap_1$ and let $ \,x_1=n_k+n_i+n_j, \quad   \,x_2= 2n_i+n_j$. \, Then \\
 \centerline{ $  Ap_2\supseteq \{  2n_i, n_i +n_j ,n_i +n_k ,n_j +n_k    \}$}.
   ($a )$.
     If the  four elements above are distinct,  we deduce that $C_3=\{n_k +n_i+n_j ,2n_i  +n_j ,3n_i,  2n_i+n_k\}$, because we must have   $  C_3\setminus\{x_1,x_2\}    \subseteq
 \{ 3n_i,  2n_i+n_k\}$ (recall that the elements of $C_3$ induce elements in $Ap_2$,  and $| C_3|\geq 4$).  \\  
  $(b_1)$.  Otherwise  $n_k+n_j=2n_i$, then $\, x_2=n_k+2n_j$ and so $2n_j\in Ap_2$. This implies \\ \centerline{$  Ap_2=\{  2n_i, n_k +n_i ,n_i +n_j ,2n_j   \}$.} In fact these elements are distinct: if not, $2n_j=n_k+n_i$ hence $x_1=3n_i=3n_j$. \, We deduce that  $  C_3\setminus\{x_1,x_2\}    \subseteq
 \{    n_k+2n_i,n_i+2n_j,3n_j \}$. In fact, \,   $n_k+2n_j=2n_i+n_j=x_2;$ \, further  $2n_i+n_k=2n_k+n_j\notin C_3$, otherwise $ 2n_k\in Ap_2 \I 2n_k=n_i+n_j$ then $x_1=3n_k=3n_i$. Hence\\ \centerline{$C_3=\{n_k +n_i+n_j=3n_i  ,2n_i  +n_j ,  n_i  +2n_j, 3n_j.\}$ } \\[2mm]
 Step {\it 5}. Now we assume that   $| Supp( x)|\leq 2\,\, \forall x\in C_3$ and  that    $ Supp\,x_1= Supp\,x_2=\{n_i,n_j\}$.\\
Let   $x_1=    n_i+2n_j$, $x_2=2n_i+ n_j$. \, Then $  Ap_2\supseteq\{2n_i,n_i+n_j,  2n_j\}$:\\[1mm] 
 $(b_2)$.   If \,$Ap_2=\{2n_i,n_i+n_j,  2n_j,n_i+n_k\}$, \,\,then   $C_3\subseteq \{3n_i,2n_i+ n_j,n_i+2n_j,3n_j,2n_i+n_k\}$. 
 \\[1mm]
 $( c \,)$.  \, If $\,Ap_2=\{2n_i,n_i+n_j,  2n_j,n_h+n_k\}$,\,\, $n_h\neq n_k$, then   $C_3=\{3n_i,2n_i+ n_j,n_i+2n_j,3n_j\}$. \\ In fact $2n_h+n_k\notin C_3$, otherwise $2n_h\in Ap_2\I 2n_h=n_i+n_j\I 2n_h+n_k=n_i+n_j+n_k $, then\\ $| Supp(2n_h+n_k)|\geq 3$, against the assumption.\\[1mm]
$(d)$  \,If\, $Ap_2=\{2n_i,n_i+n_j,  2n_j,2n_k\}$, then   $C_3\subseteq \{3n_i,2n_i+ n_j,n_i+2n_j,3n_j,3n_k\}$. \\[2mm]
Step {\it 6}. 
  Finally assume $| (Supp\,x_1 \cap Supp\ x_2)|=1$: two possible cases. 
  \\[2mm] 
  $(b_3)$\quad   $x_1=  2n_i+n_k$, $x_2=      n_i+2n_j$, then $  Ap_2=\{2n_i,n_i+n_j, n_k+n_i,2n_j\}$. In fact these four elements are distinct, otherwise
  $2n_j=n_i+n_k$ would imply $x_1= x_2$.
  Then $C_3\subseteq \{3n_i,2n_i+ n_j,n_i+2n_j,3n_j, 2n_i+n_k\}$.\\
  Note that statement $(b)$ summarizes cases $ (b_1),(b_2),(b_3)$. \\[1mm]
  $(e)$\quad  $x_1=2n_i+ n_k$, $x_2=2n_j+ n_k$ then $\{2n_i,n_i+n_k, n_j+n_k,2n_j\}= Ap_2$. In fact these four elements are distinct, otherwise
  $2n_j=n_i+n_k$ would imply $x_1=2n_i+n_k,\,\, x_2=n_i+2n_k$ and so $Supp(x_1)\cap Supp(x_2)=\{n_i,n_k\}$.
  We conclude that $ C_3= \{3n_i,2n_i+ n_k, 2n_j+n_k,3n_j \}$.   \quad$\diamond$ 

  \begin{ex}   \label{es5}{\rm We list some    semigroups  verifying   (\ref{ap24}), the second and third ones verify     also (Theorem \ref{40}).}\enu
  
  \item {\rm (\ref{ap24}.$b$) \,  Let $S=<30,33,37,73,76,77,79\rightarrow 89,91,92,94,95,98,101,108 >$ and let  $n_i=33,n_j=37,n_k=98$. Then:\\    $Ap_2= \{66,70,74,135\}=\{2n_i,n_i+n_j,2n_j,n_j+n_k\}$, $Ap_3=\{3n_i\}$, $Ap_4=\{4n_i\}$, $v=e-6$. \\ $D_2+e=\{ 103=2n_i +n_j, 107=n_i+2n_j ,\,111=3n_j  \}\subseteq C_3=\{3n_i,2n_i+n_j,n_i+2n_j,3n_j\}$ \\
  $D_3+e=\{165=5n_i, 136=3n_i+n_j, 140=2n_i+2n_j,144=n_i+3n_j,148=4n_j\},$ \\$ C_4=\{4n_i, 136,140,144,148\},$ \\ $ D_4+e=\{169=4n_i+n_j, 173=3n_i+2n_j,177=2n_i+3n_j, 181=n_i+4n_j,185=5n_j, 198=6n_i\}$.   $H_R=[ 1, 24, {\bf 25, 24, 23,} 25, 27, 29, 30\rightarrow ]$ decreases at levels  3 and 4. }
  
 \item{\rm(\ref{ap24}.$b$), with $| C_3|=5)$.  Let $S=<17,19,22,31,40,42,43,45,46, 47, 49, 52,54 >$ and let  $n_i=19,n_j=22,n_k=31$. Then:\, 
   $Ap_2= \{38,41,44,50\}$, $Ap_3=\emptyset$,  $\, v=e-4$. \\ $C_3=D_2+e =\{57=3n_i  , 60=2n_i+n_j ,\,63=n_i+2n_j  ,\,66=3 n_j , \, 69=2n_i+n_k  \}$.  Hence  $H_R$   decreases at level 2,\, $H_R=[1, {\bf13,12},14,16,17\rightarrow]$.}    
   
   \item{\rm (\ref{ap24}.$d$). Let $S=<17,22,29,37,49,64,69,70,79,82,84,89,94 >$ and let $n_i=22,n_j=37,n_k=29$. Then:  \\ $Ap_2= \{44,58,59,74\}$, $Ap_3=\emptyset$, $v=e-4$. \quad $C_3=D_2+e=\{66,81,96,  111, 87 \}=\{3n_i, 2n_i+n_j, n_i+2n_j, 3n_j, 3n_k $. \quad    $H_R=[1, {\bf13,12},14,16,17\rightarrow]$  \,\,  decreases at level 2.}
  \denu  \end{ex}  

\section{Decrease of the Hilbert function: the cases $\, v=e-3,\,\,v=e-4$.}  
\subsection{Case $\quad v=e-3$\,.}
 Assume   the embedding dimension $v$ and the multiplicity $e$ satisfy  $\,\, v=e-3$. With notation \ref{CM}, the Hilbert function  of the ring $R'=R/ t^e R$     \,is\,  $H_{R'} =[ 1,v-1,a,b,c]$\,\, with \,\,  $a+b+c=3$;  
  by Macaulay's theorem, the admissible cases are \\ 
  \centerline{$\left.\begin{array} {cllll} (i)&a=1&b=1 & c=1&(stretched \,\, case)  \\
  (ii)&a=2&b=1 & c=0 \\
 (iii)& a=3&b=0& c=0&(short \,\, case)
\end{array}
\right.$} 
 
\begin{rem}\label{gor3}  
{\rm  
In  cases $(i), (ii)$ of the above table we have \, $ | Ap_2|\leq 2  $, hence $H_R$ is non-decreasing by (\ref{pt}.5\,$c$).\\
Then  the possible {\it decreasing} examples have $H_{  R' }=[1,v-1,3]$ (short case). In any case, it is clear that no $R=k[[S]]$ with $S$ symmetric and $v=e-3$ has decreasing Hilbert function. In fact, recall that $S$ symmetric implies \,  $ \B=\{e+f  \}=Ap_m$. }
\end{rem}
\begin{thm} \label{e-3}  
With Notation \ref{def1}, assume $v=e-3$. Then the following conditions are equivalent:
\enu 
\item    $H_R$ decreases.
\item    $H_R$ decreases at level 2.
 \item  \enu
\item The Hilbert function of $R'=R/t^eR$ \, is \, $  [ 1,e-4,3]$.
\item {There exist distinct elements $n_i,  n_j \in Ap_1$  
  such that: 
  \item[]      $C_2=\{2n_i,  n_i+n_j ,   2n_j \}$   \quad  $D_2+e= C_{3}=\{3n_i, 2n_i+n_j ,n_i+2n_j ,   3n_j \}.$ } \quad 
\denu 
\denu
\end{thm}
Proof. 
If $H_R$ decreases at some level, then by (\ref{gor3}) we have $H_{R'}=[1,e-4,3]$, hence $Ap_3=\emptyset $, and so    $D_2+e= C_3$, \,  by \, (\ref{CnonAp}.1\,$b$). \\
${\it  1 \I 2 }$.  \,  If $H_R$ decreases at level $j\geq3$, then  by  \cite[Corollary 4.2]{ddm},  $| C_h|\geq h+1$ for each $2\leq h\leq j$. In particular we get $|C_3|\geq 4$, hence $|D_2|\geq 4>|C_2|=3$ and  $H_R$ decreases at level 2.  
\\[2mm]  
 ${\it  2 \I 3 }.$ If $H_R$ decreases at level 2, then by \cite[Corollary 2.4]{ddm} $| D_2|\geq 4$, hence $| C_3 |\geq 4$  and   by applying   Theorem \ref{c23}  we deduce $   D_2+e =C_3=\{3n_i, 2n_i+n_j ,n_i+2n_j ,   3n_j \}$.  
  Hence the thesis.\\[2mm]
  $  {\it  3 \I 1. }$ It is clear. \quad$\diamond$  \\[2mm]
  Note that the element of $D_2$ are an arithmetic sequence $(x_{k+1}=x_k+n_j-n_i), \, k=1,2,3 $. \\
  Under the assumption $v=e-3$, Theorem \ref{e-3} allows to prove that  $e=13$ is  the smallest multiplicty for a semigroup with decreasing Hilbert function. This bound is sharp, as shown in  Examples (\ref{tammol}) and (\ref{es13}).  
 
  \begin{prop}\label{e13} 
  If \, $v=e-3$ and the equivalent conditions of $(\ref{e-3})$ hold, then:
  \enu  
  \item There exist ten distinct elements in $S$ such that  \\[2mm]\centerline {$   \begin{array}{lll} M\setminus 2M  \supseteq  \{ e, n_i, n_j, 3n_i\! -\! e,2n_i\! +\! n_j\! -\! e, n_i\! +\! 2n_j\! -\! e, 3n_j\! -\! e\}, &
  Ap_2 =\{2n_i,  n_i+n_j ,   2n_j \}.
  \end{array} $ }
   \item If either $e$ is odd, or e=12, there exists $h\geq 2$ such that \, $(2n_i\! +\!2n_j\! -\!he) \in Ap_1$.
  \item $e=e(S)\geq 13$.
   \denu
  
\end{prop}
  Proof. {\it  1}. \,  For simplicity we denote   $n_i<n_j$ the elements $\in S$  such that $C_2=\{2n_i,n_i+n_j,2n_j\}$. By the assumption, $|D_2 | =4$, and   \ $D_2+e=\{ 3n_i, 2n_i+n_j, n_i+2n_j, 3n_j\}.$\,\, Hence\\ \centerline{ \quad
  $M\setminus 2M\supseteq  \{e, n_i, n_j, 3n_i\! -\! e,2n_i\! +\! n_j\! -\! e, n_i\! +\! 2n_j\! -\! e, 3n_j\! -\! e\}$; 
 $\quad Ap_2=
 \{ 2n_i,   n_i\! +\! n_j,  2n_j \}(=C_2)$.}
Clearly, to see that the above elements of $M\setminus 2M$ are all distinct, it's enough to verify  that $3n_i-e\neq n_j $:  this implies $
|M\setminus 2M|\geq 7 $,  $e\geq 10$. \\
 Assume   $3n_i-e=n_j$, then $ C_2=\{2n_i, 4n_i-e, 6n_i-2e\},$\, $ D_2=\{3n_i-e, 5n_i-2e, 7n_i-3e, 9n_i-4e\} $; hence   \\
  $M\setminus 2M\supseteq 
 \{e, n_i,3n_i\! -\!e,  5n_i\! -\!2e,  7n_i\! -\!3e, 9n_i\! -\!4e\}$, \, $Ap_2=
 \{ 2n_i,  4n_i\! -\!e,  6n_i\! -\!2e \}$.\, Since  $   8n_i-3e=(3n_i-e)+(5n_i-2e)\notin Ap_2, $ one has $   8n_i-4e\in M ,$ \, impossible   (
otherwise  \,  $  9n_i\! -\!4e\in 2M\cap (M\setminus 2M)$).     \\[2mm]
{\it  2}. \, Note that $n_i+(n_i+2n_j-e)=2n_i+2n_j- e\in 2M\setminus Ap_2$, then 
  $2n_i+2n_j- 2e\in S$.   It is easy to see that $2n_i+2n_j- 2e\notin  < n_i,n_j,3n_i-e, \quad 2n_i+n_j-e, \quad n_i+2n_j-e,3n_j-e  >$.\\
  Moreover for all $y\in Ap_2$ and for all $k\geq 0$, we cannot have $ 
  2n_i+2n_j- 2e=y+ke$. 
 Hence \, $2n_i+2n_j- 2e\in (Ap_1+he)\cup <e>$.\\
 Now we show that $2n_i+2n_j=\lambda e $ is impossible for $e$ odd and for $e=12$.\\
    Clearly $e\,\, odd \I \lambda$ even, and so $n_i+n_j\equiv e$, impossible since $n_i+n_j\in Ap_2$.\\
    Let $e=12$. Then $n_i+n_j=6\lambda,\,  \lambda$ odd (otherwise $n_i+n_j=\mu e$. Then $n_j\equiv -n_i-6\equiv 6+11n_i \,\, (mod\, 12)$. Clearly we cannot have \\[2mm] 
    $\overline{n_i}, \overline{n_j}\in\{\overline{3},\overline{4},\overline{6},\overline{8},\overline{9}\}$, \quad hence $\overline{(n_i,n_j)}\in\{\overline{(1,5)}  ,\overline{(5,1)},\overline{(7,11)},\overline{(11,7)}\}$. These remaining cases are impossible, because  they imply \,  $3n_i-e\equiv 3n_j-e$. Hence there exists $h\geq 2$ such that  $2n_i+2n_j-he\in (M\setminus 2M)$.
\\[2mm]     
{\it  3}.  To prove that the equivalent conditions of $(\ref{e-3})$ imply $e\geq 13$,  we  proceed   in two steps. First note that $|Ap|\geq 10$, hence $e\geq 10$.
By $(2)$, we can  assume that {$ \left\{  \begin{array}{lll} M\setminus 2M  \supseteq  \{e, n_i, n_j, 3n_i\! -\! e,2n_i\! +\! n_j\! -\! e,2n_i\! +\! 2n_j\! -\!h e, n_i\! +\! 2n_j\! -\! e, 3n_j\! -\! e\},\\
  Ap_2 =\{2n_i,  n_i+n_j ,   2n_j \}.
  \end{array}\right.$ }
\\[2mm]
{\it Step 1}: \quad $3n_i+n_j-2e\in (M\setminus 2M) $, for  $e\in \{10,11,12\} $; hence $ |$ {\it Ap\'ery} | $\geq 12$ and we cannot have $e\in \{10,11\}  $.\quad    \\ $  3n_i+n_j-e=(3n_i-e)+n_j\notin Ap_2: ($ otherwise $n_j=3n_i-e)$. Hence  $3n_i+n_j-2e\in S$.
It is easy to check that \quad $3n_i+n_j-2e\notin < n_i,n_j,3n_i-e, \quad 2n_i+n_j-e, \quad n_i+2n_j-e,3n_j-e,2n_i+2n_j-2e  >\cup (Ap_2+ke), k\geq 0$.\\  
If \, $3n_i+n_j =\lambda e, \, \lambda >4$, then  
$  {\overline  {3 n_j }}= {\overline  {3 n_j -e}}={\overline  {-9n_i} }  \quad (mod\,\,e)$. \quad  This is impossible for $e\in\{10,11,12\}$:\\
$e=10\I   {\overline {-9n_i} }=  {\overline  { n_i} }$, impossible since $n_i, 3n_j-e\in (M\setminus 2M)$.\\
$e=11\I   {\overline { -9n_i} }=  {\overline  { 2n_i} }$, impossible since $2n_i, 3n_j-e\in Ap$.\\
$e=12\I   {\overline  {-9n_i} }=  {\overline  {3 n_i} }$, impossible since $3n_i-e, 3n_j-e\in (M\setminus 2M)$.\\
 Hence there exists $k\geq 2$ such that $3n_i+ n_j-ke\in (M\setminus 2M),$. Hence   $v\geq 8$ \, for  $e\in \{10,11,12\}\quad$.  \\[2mm]
 {\it Step 2}. \,  Now   assume  $e=12$; we prove that there exists $q\geq 3$ such that  $3n_i+2n_j-qe\in (M\setminus 2M)$, hence $|${\it Ap\'ery}|$\geq 13$, therefore we cannot have $e=12$. \, 
First,  we know, by (2) and step 1,  that there exist $h,k\in{\mathbb N}$ such that $(M\setminus 2M)  \supseteq  \{e, n_i, n_j, 3n_i\! -\! e,2n_i\! +\! n_j\! -\! e, n_i\! +\! 2n_j\! -\! e,2n_i\! +\! 2n_j\! -\! he,3n_i\! +\!  n_j\! -\! ke, 3n_j\! -\! e\}$.\\ 
Let $r:=max\{h,k\}$: one can see that the  element   $3n_i+2n_j-re\in 2M\setminus (Ap_2+{\mathbb N} e)$  and that  $3n_i+2n_j-(r+1)e\notin \{e, n_i, n_j, 3n_i\! -\! e,2n_i\! +\! n_j\! -\! e, n_i\! +\! 2n_j\! -\! e,2n_i\! +\! 2n_j\! -\! he,3n_i\! +\!  n_j\! -\! ke, 3n_j\! -\! e\}$. \\
Hence it belongs to $(Ap_1+{\mathbb N} e)\cup{\mathbb N} e $. \\
 Finally it results that   $ 3n_i+2n_j\notin {\mathbb N}  e$. Otherwise,  ${\overline {3n_i-e}}= {\overline {-2n_j}}= {\overline {10 n_j}} \I  10 n_j-3 n_i= 12k\I 10 n_j=3(n_i+4k)\I   n_j=3h \I $   $h$ odd. But this implies that $ \overline {3n_i-e}=\overline {2n_j} $, impossible.   \quad  Hence there exists $q\geq r+1\geq 3 $ such that $3n_i+2n_j-qe\in (M\setminus 2M)$. This  proves that  $e\geq 13$. \quad$\diamond$ 
  \begin{coro}\label{coro13}    If \, $v=e-3$ \, and \, either $|Ap_2| \leq 2,$\,\, $or \,\,( |  Ap_2 |=3$,\quad $e\leq 12)$ \quad   the  Hilbert function is non-decreasing. \end{coro} 
  In the next example, for $v=e-3$,  we exhibit a technique of  computation which allows at first to verify that $H_R$ decreasing $\I e\geq 13$ and   further     to give a complete  description of the semigroups $S$ with $e=13=v+3$ and decreasing Hilbert function. Similar tables could be used for $e>13$ and also if $v=e-r$, $r\geq 4$  to find semigroups with $H_R$ decreasing.
    
  \begin{ex}\label{ex13} 
   {\rm Let $v=e-3$ and $H_R$ decreasing; $e\geq 10$      (\ref{e13})     and there exist $n_i,n_j\in Ap_1,$ distinct elements such that    \\ \centerline  {$  \begin{array}{lllllllll} Ap_1  \supseteq  \{e, n_i, n_j, 3n_i\! -\! e,2n_i\! +\! n_j\! -\! e, n_i\! +\! 2n_j\! -\! e, 3n_j\! -\! e\},&&Ap_2 =\{2n_i,  n_i+n_j ,   2n_j \}   \end{array}  $. }       
  Let   $GCD(e,n_i)=1$; \,  the following table is   useful   to find the  pairs $(n_i,n_j)$ \lq\lq compatible" with the   assumption on the Ap\'ery  set.   In the table we fix $n_j\equiv h \,n_i (mod\, e)$: in the columns we indicate  the classes of  elements of the Ap\'ery  set $(mod\, n_i)$:\,\, in each row  we must have    distinct values  $(mod\, e)$. \\ 
  Under our conditions,  we consider $ 4\leq h\leq e-1$  and\quad  $10\leq e\leq 13$). Clearly,  for $e>10$ some element in $ \{2n_i+2n_j-\lambda_1 e, 3n_i+ n_j-\lambda_2 e, n_i+3n_j-\lambda_3 e,3n_i+2n_j-\lambda_4 e, 2n_i+3n_j-\lambda_5 e\}$ must belong to $Ap_1$, for this reason we add 5 columns useful to complete the Ap\'ery  set in cases $11\leq e\leq 13$.  \\[2mm] 
\centerline{  $ \begin{array}{ccccccccccccccccccccccccccc}
 n_i\!\!&\!  2n_i\!\!&\!3n_i \!\!&\! n_j\!\!&\!n_i\! +\! n_j\!\!&\! 2n_i\! +\! n_j \!\!&\!  \! 2n_j\!\!&\!n_i\! +\! 2n_j    \!\!&\! 3n_j  \!\!&\!   2n_i\! +\! 2n_j  \!\!&\!3n_i\! +\!n_j\!\!&\!n_i\! +\!3n_j \!\!&\!3n_i\! +\! 2n_j\!\!&\!2n_i\! +\! 3n_j \\
 1\!\!&\!    2 \!\!&\!3  \!\!&\! h  \!\!&\!(h\! +\!1)  \!\!&\! (h\! +\!2)  \!\!&\! 2h  \!\!&\! (2h\! +\!1) \!\!&\! 
 3h   \!\!&\!    (2h\! +\!2)  \!\!&\!(h\! +\!3)  \!\!&\!(3h\! +\!1)\!\!&\!(2h+3)\!\!&\!(3h+2)    \\
1\!\!&\!       2\!\!&\!3\!\!&\! 4\!\!&\! 5\!\!&\! 6\!\!&\! 8\!\!&\!9\!\!&\!12  \!\!&\!   10  \!\!&\!7\!\!&\! 13 \!\!&\!11\!\!&\!14\!\!&\!ok\\
1\!\!&\!      {\bf2}\!\!&\!3\!\!&\! 5\!\!&\! 6\!\!&\! 7\!\!&\!10\!\!&\!11\!\!&\!{\bf 15}  \!\!&\!  12   \!\!&\!8\!\!&\!16\!\!&\!\!\!&\!\!\!&\!no\\
1\!\!&\!   2\!\!&\!3\!\!&\! 6\!\!&\! 7\!\!&\! 8\!\!&\!12\!\!&\!{\bf13}\!\!&\!18  \!\!&\!   { 14}  \!\!&\!9\!\!&\!19\!\!&\!\!\!&\!\!\!&\!no\\
1\!\!&\!    {\bf2}\!\!&\!3\!\!&\! 7\!\!&\! 8\!\!&\! 9\!\!&\!14\!\!&\!{\bf15}\!\!&\!21  \!\!&\!   16   \!\!&\!10\!\!&\!22\!\!&\!\!\!&\!\!\!&\!no\\
 1\!\!&\!     2\!\!&\!{\bf 3 }\!\!&\! 8\!\!&\! 9\!\!&\!10\!\!&\!{\bf16 }\!\!&\!17\!\!&\!24  \!\!&\!   18  \!\!&\!11\!\!&\!25\!\!&\!\!\!&\!\!\!&\!no\\
1\!\!&\!    2\!\!&\!3\!\!&\! 9\!\!&\!10\!\!&\!11\!\!&\!18\!\!&\!19\!\!&\!{\bf 27}  \!\!&\!   20  \!\!&\!12\!\!&\!28\!\!&\!\!\!&\!\!\!&\!no\\
 1\!\!&\!     2\!\!&\!3\!\!&\!10\!\!&\!11\!\!&\!12\!\!&\!20\!\!&\!21\!\!&\!30  \!\!&\!   22  \!\!&\!13\!\!&\!5\!\!&\!23\!\!&\!32\!\!&\!ok\\
 1\!\!&\!     2\!\!&\!3\!\!&\!11\!\!&\!12\!\!&\!{\bf13 }\!\!&\!22\!\!&\!23\!\!&\!33  \!\!&\!   24  \!\!&\!14\!\!&\!34\!\!&\!\!\!&\!\!\!&\!no\\
 1\!\!&\!     2\!\!&\!3\!\!&\!12\!\!&\!{\bf13 }\!\!&\!14\!\!&\!24\!\!&\!25\!\!&\!36  \!\!&\!   26  \!\!&\!15\!\!&\!37\!\!&\!\!\!&\!\!\!&\!no\\
  \end{array}
  $}
  \\[2mm]
 The   table shows that for each $e$,  only few cases  with $H_R$ decreasing  are  "admissible". Moreover with some other check in cases  $e\in \{10,11,12\}$,  one can confirm that  $  e\geq 13$, as proved in  (\ref{e13}.3);\\[2mm]  in case $e=10$ the only remaining case is  $(n_i=2, n_j=5)$, impossible since this would imply $2n_j\equiv 0$;\\[2mm]
  in case $e=12$ the   cases to be veryfied are  $(n_i,n_j)\in\{( 2, n_j=3)(2,  9),(4,-),(6,-), (8,-)\}$, which are clearly incompatible with the assumptions \, (for $(2,9): 2n_j\equiv 3n_i$).\\[2mm]
  For $e=13$, the possible cases are $n_j\equiv 4n_i,$ \, or   \, $ n_j\equiv 10n_i$ (in the table the pairs of distinct elements with the same class $mod \, n_i$  for $e=13$ are written in bold). \\[2mm]  
   By means of these computations we deduce the structure of such semigroup rings   with $H_R$ decreasing:}
 \end{ex}
     \begin{prop} \label{exe13} 
 With Notation \ref{def1}, let $R=k[[S]]$, with  multiplicity $e=  13$ and $v=10(=e-3)$. Further let $1\leq p\leq 12:$  there exists  $R_p$ satisfying the equivalent conditions of Theorem \ref{e13} if and only if there exist $(k,k',\alpha,\beta,\gamma)_p$, \, $k,\alpha,\beta,\gamma \in{\mathbb N},\, k'\in{\mathbb Z}$,  such that  semigroup $S_p$ has   the following minimal set of generators.  
 \enu \item[]    
 $S_p=<  e, n_i, n_j, 3n_i\! -\! e,2n_i\! +\! n_j\! -\! e, n_i\! +\! 2n_j\! -\! e, 3n_j\! -\! e,   2n_i\! +\! 2n_j\! -\! \alpha e,    3n_i\! +\!  n_j\! -\! \beta e,3n_i\! +\!  2n_j\! -\! \gamma e> $ 
 \item[]  with\,\, $n_i=k e+p,\,\, k \geq 1,\,\,    n_j=k'e+ 4p,\quad    -2\leq k'\leq 4k-2, \quad 4k'>3k-p,\quad\alpha<\gamma,\quad \beta<\gamma. $
 \denu
   \end{prop}
   Proof.  By (\ref{ex13}), for $e=13$, the possible cases are $n_j\equiv 4n_i,$ \, or   \, $ n_j\equiv 10n_i$: these conditions are symmetric, $n_i\equiv 4n_j \II n_j\equiv 10 n_i$, hence there is only one class of semigroups with decreasing Hilbert function.
\quad  Further :\\
    -   $   n_j=k'e+4p=4n_i-re \I r\geq 2$ (otherwise $r=1\I n_j=(3n_i-e)+n_i\in 2M)$.  \\
 -   $3n_i\equiv 4n_j \I 3n_i-e< n_j+(3n_j-e)\I 3n_i<4n_j, i.e., \,\, 3k<4k'+p$. 
    \\ The remaining 3 generators must be equivalent to $  2n_i\! +\! 2n_j,3n_i\! +\!  n_j, 3n_i\! +\!  2n_j\,(mod\,e).$  Let $2n_i+2n_j-\alpha e\in Ap_1$, then $\alpha\geq 2$ ($\alpha=1\I 2n_i+2n_j- e=n_j+(2n_i+n_j-e)\notin Ap_1$); now note that $3n_i\! +\!  2n_j-\gamma e =(2n_i+2n_j-\alpha e)+ n_i +(\alpha-\gamma) e\notin M+e\I \gamma>\alpha$; analogously, $3n_i\! +\!  2n_j-\beta e=(3n_i+ n_j-\beta e)+n_j+(\beta-\gamma) \I \gamma>\beta$.\quad$\diamond$\\[2mm]
{   We exhibit, in case  $e=13=v+3$, for each $1\leq p\leq 12$ an example of semigroup $S_p$ as in (\ref{exe13}).
      \begin{ex}\label{es13}  
       {\rm  With the notation of (\ref{exe13}) above, let $k=1$, \,\, $1\leq p\leq 12$;   the following semigroups $S_p$ with  $k'$ minimal ($-2\leq k'\leq 1$) have decreasing Hilbert function  ($S_6$ is    Example
     \ref{tammol} ).
\\[2mm]     $\begin{array}{llccllll}
  ({\overline n_i},{\overline n_j})=({\overline 1},{\overline 4)}&:\,\, S_1=<13,14,17,29, 32,33, 35,36, 37,38>& k'= 1&  c=26\\ 
 ({\overline n_i},{\overline n_j})=({\overline 2},{\overline 8)} &:\,\,  S_2=<13,15,21,32,38,40,44,46,48,50>  & k'=1    &c=38\\ 
   ({\overline n_i},{\overline n_j})=({\overline 3},{\overline {12}}), &:\,\,  S_3=<13,16,25,35, 44,47,53,56,59,62> &  k'=1   &c=50\\
   ({\overline n_i},{\overline n_j})=({\overline 4},{\overline 3)} &:\,\,  S_4=<13,17,16,35, 36,37 ,38 , 40,41,44>& k'=0  &c=32\\

   ({\overline n_i},{\overline n_j})=({\overline 5},{\overline 7)} &:\,\,  S_5=<13,18,20,41,43, 45,47,48,50,55> &  k'=0  &c= 43\\
   ({\overline n_i},{\overline n_j})=({\overline 6},{\overline {11}}) &:\,\,  S_6=<13,19,24,44,49,54,55,59,60,66> &  k'=0  &c=54 \\   ({\overline n_i},{\overline n_j})=({\overline 7},{\overline 2)} &:\,\,  S_7=<13,20,28,47,55,62,63,70,71,77> &  k'=0  &c=65\\ 
   ({\overline n_i},{\overline n_j})=({\overline 8},{\overline 6)} &:\,\,  S_8=<13,21, 19, 44,46,48,50,54,56,62 > &   k'=-1  &c=50\\
      ({\overline n_i},{\overline n_j})=({\overline 9},{\overline {10}}) &:\,\, S_9=<13,22,23,53,54,55,56,63,64,73> &  k'=-1  &c=61\\
 ({\overline n_i},{\overline n_j})=({\overline 10},{\overline 1)} &:\,\,   S_{10}=<13,23,27,56,60,64,68,70,74,84> &  k'=-1  &c=72\\
 ({\overline n_i},{\overline n_j})=({\overline 11},{\overline 5)} &:\,\,   S_{11}=<13,24,31,59,66,73,77,80,84,95 > &  k'=-1  &c=83\\
 ({\overline n_i},{\overline n_j})=({\overline 12},{\overline 9)} &:\,\,   S_{12}=<13,25,22,53, 56,59,62,68, 71,80  > &  k'=-2  &c=68
\end{array}$.   }\end{ex}}
 \subsection{Case $\quad v=e-4$\,.}
With notation \ref{CM} and  \ref{def1},    the Hilbert function of   $R'=R/t^eR$,      is  {   $H_{R'}(z)\,= [1,v-1,a,b,c,d]$ \,\, with \,$a+b+c+d=4.$  
\begin{(*)} \label{star2} 
{\rm By Macaulay's Theorem the admissible   $H_{\!R'}$,  are\\[2mm] \centerline{$\left.\begin{array} {lllrr}  &[1,v-1,4]  \\ 
 &[1,v-1,3,1]  \\
  &[1,v-1,2,2] \\
 &[1,v-1,2,1,1] \\
 &[1,v-1,1,1,1,1] &(stretched)
\end{array}
\right.$}\\[2mm]
When $|  Ap_2 |\leq 2 $ we know by (\ref{pt}.5\,$c$) that $H_R$ is non-decreasing.
Hence we consider the first two cases.}\end{(*)}
  \begin{thm} \label{40} 
   With Notation \ref{def1}, assume $v=e-4,\, |  Ap_2 |=4$, $  Ap_3=
\emptyset$. The  following conditions are equivalent:
\enu 
 \item $H_R$ decreases at level 2.
\item There exist $n_i,n_j,n_k\in Ap_1$, distinct elements, such that  
\item[] either $Ap_2=   \{ 
 2n_i,
 n_i+n_j, 
2n_j, 
n_i+n_k \}, \quad  
  C_3=D_2+e= \left\{   3n_i ,
  2n_i+n_j,
  n_i+2n_j,
  3n_j,
  2n_i+n_k )
  \right\}$
  \item[] or \qquad    $Ap_2=   \{ 
 2n_i,
 n_i+n_k, 
2n_j, 
  2n_k)\}, \quad  
  C_3=D_2+e= \left\{   3n_i ,
  2n_i+n_j,
  n_i+2n_j,
  3n_j,
 3n_k)
  \right\}$

\denu
\end{thm}
 Proof. By  (\ref{CnonAp}.1\,$b$) we have $\, C_3=D_2+e$, then  $H_R$ decreases at level 2 $\II$ $|D_2 |\geq 5$, i.e. $ | C_3 |\geq 5$; now apply (\ref{ap24}).\quad$\diamond$

 \begin{thm} \label{31}    With Notation \ref{def1}, assume    $v=e-4, |  Ap_2 |=3$, $   | Ap_3 |=1$. The  following conditions are equivalent:
\enu 
 \item $H_R$ decreases.
  \item $H_R$ decreases at level $h\leq 3$.
 
\item  There exist $n_i,n_j \in Ap_1$, distinct elements, such that \\  $ Ap_2 =  \left\{\begin{array}{lllll}
 2n_i \\
 n_i+n_j  \\
2n_j\\
 \end{array}\right.\!\!,\,\,C_3= \left\{\begin{array}{lllll}
  3n_i\in Ap_3 \\
  2n_i+n_j\\
  n_i+2n_j\\
  3n_j
 \end{array}\right., \,  \, D_h+e= \left\{\begin{array}{lllll}
   4 n_i  \\
 hn_i+n_j \\
 (h-1) n_i+2n_j \\
 \cdots \\
  (h+1)n_j 
 \end{array}\right.,\quad (h=2\,or\, h=3).$   
 
\denu
\end{thm}
 Proof. ${\it 1\I2}$, and ${\it1\I 3}$ follow by (\ref{apj1}.4), with $d=3$.\\
   ${\it3\I2\I1}$ are obvious.\quad$\diamond$

\begin{coro}\label{corogor} 
 If $\, R=k[[S]]$ is Gorenstein with $v\geq e-4$, \, then $H_R$ is non decreasing.
\end{coro}
Proof.   For $v\geq e-2$, see \cite{e1}\, \cite{e2}. For $v=e-3$ see   \ref{gor3}. In case $v=e-4$ the result follows by (\ref{star2}), (\ref {pt}.5$c$),  the properties of symmetric semigroups and 
  (\ref{apj1}.4$b$). \quad$\diamond$\\[2mm]
  
\noindent The authors are grateful to Professor Francesco Odetti Dime-University of Genova for his essential informatic support.

\end{document}